\documentclass[11pt]{article}
\usepackage{epic,latexsym,amssymb}
\usepackage{color}
\usepackage{tikz}

\newtheorem{thm}{Theorem}
\newtheorem{lem}[thm]{Lemma}

\newtheorem{cor}{Corollary}
\newtheorem{ob}{Observation}
\newtheorem{prop}{Proposition}

\newcommand{\qed}{$\Box$}

\newcommand{\smallqed}{{\tiny ($\Box$)}}

\newcommand{\cC}{{\cal C}}

\newcommand{\cG}{{\cal G}}
\newcommand{\cH}{{\cal H}}

\newcommand{\oc}{{\rm oc}}

\newenvironment{unnumbered}[1]{\trivlist \item [\hskip \labelsep {\bf
#1}]\ignorespaces\it}{\endtrivlist}

\def\vertex(#1){\put(#1){\circle*{2}}}
\def\vertexo(#1){\put(#1){\circle{2}}}
\def\vert(#1){\put(#1){\circle*{1.5}}}
\def\verto(#1){\put(#1){\circle{1.5}}}
\def\lab(#1)#2{\put(#1){\makebox(0,0)[c]{#2}}}
\newcommand{\proof}{\noindent\textbf{Proof. }}
\newcommand{\3}{ \vspace{0.3cm} }
\newcommand{\2}{ \vspace{0.2cm} }
\newcommand{\1}{ \vspace{0.1cm} }

\let\oldenumerate\enumerate
\renewcommand{\enumerate}{
  \oldenumerate
  \setlength{\itemsep}{1pt}
  \setlength{\parskip}{0pt}
  \setlength{\parsep}{0pt}
}

\begin{document}

\title{Tight lower bounds on the matching number \\ in a graph with given maximum degree}

\author{Michael A. Henning${}^{1,}$\thanks{Research
supported in part by the South African
National Research Foundation and the University of Johannesburg} \, and Anders Yeo${}^{1,2}$\\
\\
${}^1$Department of Pure and Applied Mathematics\\
University of Johannesburg \\
Auckland Park, 2006 South Africa \\
Email:  mahenning@uj.ac.za \\
\\
${}^2$Engineering Systems and Design \\
Singapore University of Technology and Design \\
8 Somapah Road, 487372, Singapore \\
Email: andersyeo@gmail.com
}

\date{}
\maketitle

\begin{abstract}
Let $k \geq 3$. We prove the following three bounds for the matching number, $\alpha'(G)$, of a graph, $G$, of order $n$ size $m$ and maximum degree at most $k$.

\begin{itemize}
\item If $k$ is odd, then $\alpha'(G) \ge \left( \frac{k-1}{k(k^2 - 3)} \right) n \, + \, \left( \frac{k^2 - k - 2}{k(k^2 - 3)} \right) m \, - \, \frac{k-1}{k(k^2 - 3)}$.
\item If $k$ is even, then  $\alpha'(G)  \ge  \frac{n}{k(k+1)} \, + \, \frac{m}{k+1} - \frac{1}{k}$.
\item If $k$ is even, then  $\alpha'(G)  \ge  \left( \frac{k+2}{k^2+k+2} \right) m  \, - \,  \left( \frac{k-2}{k^2+k+2} \right) n \, - \frac{k+2}{k^2+k+2}$.
\end{itemize}

In this paper we actually prove a slight strengthening of the above for which the bounds are tight for essentially all densities of graphs.

The above three bounds are in fact powerful enough to give a complete description of the set $L_k$ of pairs $(\gamma,\beta)$ of real numbers with the following property.
There exists a constant $K$  such that $\alpha'(G) \geq \gamma n + \beta m - K$ for every connected graph $G$ with maximum degree at most~$k$,
where $n$ and $m$ denote the number of vertices and the number of edges, respectively, in $G$.
We show that $L_k$ is a convex set.
Further, if $k$ is odd, then $L_k$ is the intersection of two closed half-spaces, and there is exactly one extreme point of $L_k$, while if $k$ is even, then $L_k$ is the intersection of three closed half-spaces, and there are precisely two extreme points of $L_k$.
\end{abstract}

{\small \textbf{Keywords:} Matching number, maximum degree, convex set} \\
\indent {\small \textbf{AMS subject classification:} 05C65}

\newpage
\section{Introduction}
\label{Intro}

Two edges in a graph $G$ are {\em independent} if they are not adjacent in $G$. A set of pairwise independent edges of $G$ is called a {\em matching} in $G$, while a matching of maximum cardinality is a {\em maximum matching}. The number of edges in a maximum matching of $G$ is called the {\em matching number} of $G$ which we denote by $\alpha'(G)$. %
Matchings in graphs are extensively studied in the literature (see, for example, the classical book on matchings my Lov\'{a}sz and Plummer~\cite{LoPl86}, and the excellent
survey articles by Plummer~\cite{Pl03} and Pulleyblank~\cite{Pu95}).

For $k \ge 3$, let $L_k$ be the set of all pairs $(\gamma,\beta)$ of real numbers for which there exists a constant $K$  such that
\[
\alpha'(G) \ge \gamma n + \beta m - K
\]
holds for every connected graph $G$ with maximum degree at most~$k$, where $n$ and $m$ denote the number of vertices and the number of edges, respectively, in $G$. Our main result is to give a complete description of the set $L_k$. For this purpose, let $\ell_1$, $\ell_2$, $\ell_3$ and $\ell_4$ be the following four closed half-spaces over the reals $\gamma$ and $\beta$:
\[
\begin{array}{lrcccccr}
\ell_1: & \beta & \le & -\gamma & + & \frac{1}{k} \2 \\
\ell_2: & \beta & \le & - \left(\frac{2}{k} \right) \gamma & + & \frac{k^3 - k^2 - 2}{k^2(k^2 - 3)} \2 \\
\ell_3: & \beta & \le & - \left( \frac{2}{k} \right) \gamma & + & \frac{k^2 + 4}{k(k^2 + k + 2)} \2 \\
\ell_4: & \beta & \le & -\left(\frac{2k^2 }{k^3-k+2} \right) \gamma & + & \frac{k^2 - k + 2}{k^3-k+2}
\end{array}
\]

We are now in a position to state our main result.

\begin{unnumbered}{Theorem~A}
For $k \ge 3$, the set $L_k$ is a convex set. Further, the following holds.
\\[-20pt]
\begin{enumerate}
\item If $k \ge 3$ is odd, then $L_k$ is the intersection of the two closed half-spaces $\ell_1$ and $\ell_2$, and there is exactly one extreme point of $L_k$, namely
\[
\left( \frac{k-1}{k(k^2 - 3)},\frac{k^2 - k - 2}{k(k^2 - 3)} \right).
\]
\item If $k \ge 4$ is even, then $L_k$ is the intersection of the three closed half-spaces $\ell_1$, $\ell_3$ and $\ell_4$, and there are precisely two extreme points of $L_k$, namely
\[
\left( \frac{1}{k(k+1)}, \frac{1}{k+1} \right)\hspace*{0.5cm} \mbox{and}  \hspace*{0.5cm} \left( - \, \frac{k - 2}{k^2 + k + 2}, \frac{k+2}{k^2 + k + 2} \right).
\]
\end{enumerate}
\end{unnumbered}

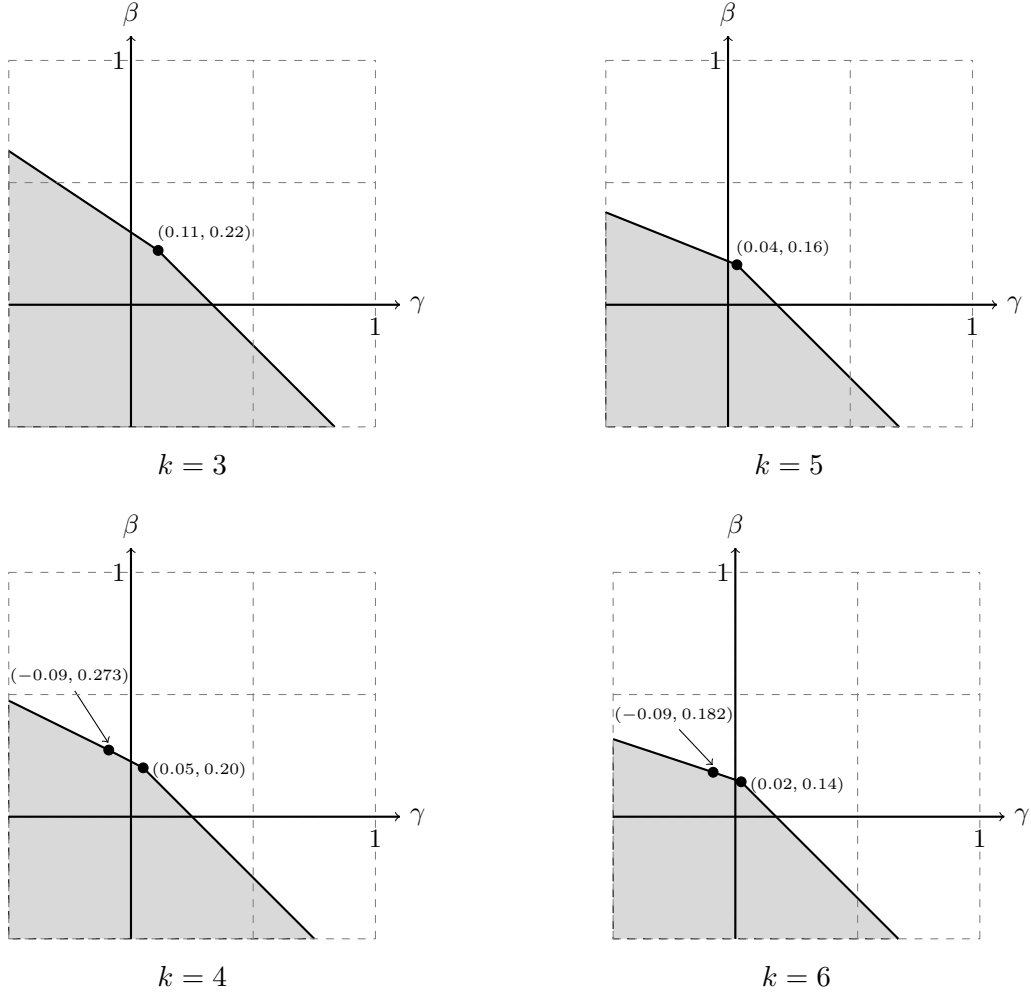
\begin{figure}[htb]
\begin{center}
\begin{tikzpicture}[scale=3.25]
\draw[fill=gray!30, dashed] (0.8333333333333333,-0.5) -- (0.1111111111111111,0.2222222222222222) --  (-0.5,0.6296296296296297) -- (-0.5,0.6296296296296297) -- (-0.5,-0.5) -- (0.8333333333333333,-0.5);
\draw[thick] (0.8333333333333333,-0.5) -- (0.1111111111111111,0.2222222222222222) --  (-0.5,0.6296296296296297);
\draw (0.30111,0.29222) node {{\tiny $(0.11,0.22)$}};
\draw[circle, fill=black]  (0.1111111111111111,0.2222222222222222) circle (0.02cm);
\draw[thin,color=gray,step=0.5cm,dashed] (-0.5,-0.5) grid (1,1);
\draw[->] (0.0,-0.5) -- (0.0,1.1) node[above] {{\small $\beta$}};
\draw[->] (-0.5,0.0) -- (1.1,0.0) node[right] {{\small $\gamma$}};
\draw[thick] (0.0,-0.5) -- (0.0,1.1);
\draw[thick] (-0.5,0.0) -- (1.1,0.0);
\draw (-0.05,1) node {{\small $1$}};
\draw (1,-0.09) node {{\small $1$}};
\draw (0.25,-0.65) node {$k=3$};
\end{tikzpicture}
 %
 \hspace*{2cm}
 %
\begin{tikzpicture}[scale=3.25]
\draw[fill=gray!30, dashed] (0.7,-0.5) -- (0.03636363636363636,0.16363636363636364) --  (-0.5,0.3781818181818182) -- (-0.5,0.3781818181818182) -- (-0.5,-0.5) -- (0.7,-0.5);
\draw[thick] (0.7,-0.5) -- (0.03636363636363636,0.16363636363636364) --  (-0.5,0.3781818181818182);
\draw (0.22636,0.23364) node {{\tiny $(0.04,0.16)$}};
\draw[circle, fill=black]  (0.03636363636363636,0.16363636363636364) circle (0.02cm);
\draw[thin,color=gray,step=0.5cm,dashed] (-0.5,-0.5) grid (1,1);
\draw[->] (0.0,-0.5) -- (0.0,1.1) node[above] {{\small $\beta$}};
\draw[->] (-0.5,0.0) -- (1.1,0.0) node[right] {{\small $\gamma$}};
\draw[thick] (0.0,-0.5) -- (0.0,1.1);
\draw[thick] (-0.5,0.0) -- (1.1,0.0);
\draw (-0.05,1) node {{\small $1$}};
\draw (1,-0.09) node {{\small $1$}};
\draw (0.25,-0.65) node {$k=5$};
\end{tikzpicture}
\vskip 0.25cm
\begin{tikzpicture}[scale=3.25]
\draw[fill=gray!30, dashed] (0.75,-0.5) -- (0.05,0.2) -- (-0.09090909090909091,0.2727272727272727) -- (-0.5,0.47500000000000003) -- (-0.5,0.47500000000000003) -- (-0.5,-0.5) -- (0.75,-0.5);
\draw[thick] (0.75,-0.5) -- (0.05,0.2) -- (-0.09090909090909091,0.2727272727272727) -- (-0.5,0.47500000000000003);
\draw[circle, fill=black]  (-0.09090909090909091,0.2727272727272727) circle (0.02cm);
\draw (0.28000,0.19000) node {{\tiny $(0.05,0.20)$}};
\draw (-0.25000,0.57500) node {{\tiny $(-0.09,0.273)$}};
\draw[thin,->] (-0.23,0.5150000000000001) -- (-0.09290909090909091,0.30272727272727273);
\draw[circle, fill=black]  (0.05,0.2) circle (0.02cm);
\draw[thin,color=gray,step=0.5cm,dashed] (-0.5,-0.5) grid (1,1);
\draw[->] (0.0,-0.5) -- (0.0,1.1) node[above] {{\small $\beta$}};
\draw[->] (-0.5,0.0) -- (1.1,0.0) node[right] {{\small $\gamma$}};
\draw[thick] (0.0,-0.5) -- (0.0,1.1);
\draw[thick] (-0.5,0.0) -- (1.1,0.0);
\draw (-0.05,1) node {{\small $1$}};
\draw (1,-0.09) node {{\small $1$}};
\draw (0.25,-0.65) node {$k=4$};
\end{tikzpicture}
 %
 \hspace*{2cm}
 %
\begin{tikzpicture}[scale=3.25]
\draw[fill=gray!30, dashed] (0.6666666666666666,-0.5) -- (0.023809523809523808,0.14285714285714285) -- (-0.09090909090909091,0.18181818181818182) -- (-0.5,0.31746031746031744) -- (-0.5,0.31746031746031744) -- (-0.5,-0.5) -- (0.6666666666666666,-0.5);
\draw[thick] (0.6666666666666666,-0.5) -- (0.023809523809523808,0.14285714285714285) -- (-0.09090909090909091,0.18181818181818182) -- (-0.5,0.31746031746031744);
\draw[circle, fill=black]  (-0.09090909090909091,0.18181818181818182) circle (0.02cm);
\draw (0.25381,0.13) node {{\tiny $(0.02,0.14)$}};
\draw (-0.25000,0.41746) node {{\tiny $(-0.09,0.182)$}};
\draw[thin,->] (-0.23,0.3574603174603174) -- (-0.09290909090909091,0.21181818181818182);
\draw[circle, fill=black]  (0.023809523809523808,0.14285714285714285) circle (0.02cm);
\draw[thin,color=gray,step=0.5cm,dashed] (-0.5,-0.5) grid (1,1);
\draw[->] (0.0,-0.5) -- (0.0,1.1) node[above] {{\small $\beta$}};
\draw[->] (-0.5,0.0) -- (1.1,0.0) node[right] {{\small $\gamma$}};
\draw[thick] (0.0,-0.5) -- (0.0,1.1);
\draw[thick] (-0.5,0.0) -- (1.1,0.0);
\draw (-0.05,1) node {{\small $1$}};
\draw (1,-0.09) node {{\small $1$}};
\draw (0.25,-0.65) node {$k=6$};
\end{tikzpicture}
\end{center}
\vskip -0.5 cm
\caption{The convex set $L_k$ for small $k$}
\label{fig:Lk}
\end{figure}

\medskip
Theorem~A is illustrated in Figure~\ref{fig:Lk} for small values of $k$, namely $k \in \{3,4,5,6\}$, where the convex set $L_k$ corresponds to the grey area in the pictures.

In order to prove Theorem~A, we shall prove the following two key results on the matching number.

\begin{unnumbered}{Theorem~B}
If $k \ge 3$ is an odd integer and $G$ is a connected graph of order $n$, size $m$, and with maximum degree $\Delta(G) \le k$, then
\[
\alpha'(G) \ge \left( \frac{k-1}{k(k^2 - 3)} \right) n \, + \, \left( \frac{k^2 - k - 2}{k(k^2 - 3)} \right) m \, - \, \frac{k-1}{k(k^2 - 3)}.
\]
Further, this lower bound is achieved for infinitely many trees, and for infinitely many $k$-regular graphs.
\end{unnumbered}

In fact, Theorem~B is tight for essentially all possible densities of connected graphs with maximum degree~$k$.

\newpage
\begin{unnumbered}{Theorem~C}
If $k \ge 4$ is an even integer and $G$ is a connected graph of order $n$, size $m$ and maximum degree $\Delta(G) \le k$, then the following holds.
\[
\begin{array}{llcl}
{\rm (a)} & \alpha'(G) & \ge & \displaystyle{ \frac{n}{k(k+1)} \, + \, \frac{m}{k+1} - \frac{1}{k} } \3 \\
{\rm (b)} & \alpha'(G) & \ge & \displaystyle{ \left( \frac{k+2}{k^2+k+2} \right) m  \, - \,  \left( \frac{k-2}{k^2+k+2} \right) n \, - \frac{k+2}{k^2+k+2} }
\end{array}
\]
\end{unnumbered}

We will later see how to slightly improve the bounds in Theorem~C, such that they are also achieved for infinitely many trees
and for infinitely many $k$-regular graphs and essentially for all possible densities in between.

\section{Known Matching Results}



We shall need the following theorem of Berge~\cite{Berge} about the matching number of a graph, which is sometimes referred to as the Tutte-Berge formulation for the matching number.

\begin{thm}{\rm (Tutte-Berge Formula)}
For every graph $G$,
\[
\alpha'(G) = \min_{X \subseteq V(G)} \frac{1}{2} \left( |V(G)| + |X| - \oc(G - X)  \right).
\]
\label{Berge}
\end{thm}

An elegant proof of Theorem~\ref{Berge} was given by West~\cite{West11}. We remark that as a consequence of the Tutte-Berge Formula, it is well-known that if $X$ is a proper subset of vertices of $G$ such that $(|V(G)| + |X| - \oc(G - X))/2$ is minimum, then every odd component $C$ of $G$ contains an almost perfect matching; that is, $\alpha'(C) = (|V(C)| - 1)/2$.

The following results from \cite{HeYe07} establishes a tight lower bound on the matching number of a regular graph.

\begin{thm}{\rm (\cite{HeYe07})}
For $k \ge 2$ even, if $G$ is a connected $k$-regular graph of order~$n$, then
\[
\alpha'(G) \ge \min \left\{  \left( \frac{k^2 + 4}{k^2 + k + 2}
\right) \times \frac{n}{2}, \frac{n-1}{2} \right\},
\]
and this bound is tight. \label{KregEVEN}
\end{thm}

\begin{thm}{\rm (\cite{HeYe07})}
For $k \ge 3$ odd, if $G$ is a connected $k$-regular graph of
order~$n$, then
\[
\alpha'(G) \ge \frac{(k^3-k^2-2) \, n - 2k +
2}{2(k^3-3k)},
\]
and this bound is tight. \label{KregODD}
\end{thm}

For small values of $k \ge 3$, the results of Theorem~\ref{KregEVEN} and~\ref{KregODD} are summarized in Table~1.

{\small

\[
\begin{array}{|c||c|c|c|c|c|c|} \hline
\multicolumn{7}{|c|}{\mbox{$G$ is a connected $k$-regular graph}} \\ \hline \hline
k & 3 & 4 & 5 & 6 & 7 & 8 \\ \hline %
& & & & & &  \\
\alpha'(G) \ge & \displaystyle{ \frac{4n-1}{9} } & \displaystyle{
\min \left\{ \frac{5n}{11}, \frac{n-1}{2} \right\} } &
\displaystyle{ \frac{49n-4}{110} } & \displaystyle{  \min \left\{
\frac{5n}{11}, \frac{n-1}{2} \right\} } & \displaystyle{
\frac{73n-3}{161} } & \displaystyle{ \min \left\{ \frac{17n}{37},
\frac{n-1}{2} \right\} } \\
& & & & & &  \\\hline
\end{array}
\]
\begin{center}
\textbf{Table~1.} Tight lower bounds on the matching number of a $k$-regular, connected graph.
\end{center}

}

\section{Main Results}

Our main result establishes tight lower bounds on the matching number of a graph in terms of its maximum degree, order, size, and number of components. Our first result is the following result, a  proof of which is given in Section~\ref{S:match_thm1}.

\begin{thm}\label{match_thm}
Let $k \ge 3$ be an integer and let $G$ be a graph with $c$ components and of order $n$ and size $m$ and maximum degree $\Delta(G) \le k$.
If no component of $G$ is $k$-regular, then
\[
\alpha'(G) \ge \left\{
\begin{array}{lc}
\displaystyle{ \left( \frac{1}{k(k+1)} \right) (n - c) \, + \, \left( \frac{1}{k+1} \right) m } & \mbox{if $k$ is even} \3 \\
\displaystyle{ \left( \frac{k-1}{k(k^2 - 3)} \right) (n - c) \, + \, \left( \frac{k^2 - k - 2}{k(k^2 - 3)} \right) m } & \mbox{if $k$ is odd}.
\end{array}
\right.
\]
\end{thm}

For small values of $k \ge 3$, the results of Theorem~\ref{match_thm} are summarized in Table~2.

{\small

\[
\begin{array}{|c|rcrcrcrcc|cc|} \hline
k & \multicolumn{8}{c}{\mbox{Exact bound} } & & & \mbox{Approximate equivalent to} \\ \hline
3 &  9 \, \alpha'(G)  & \ge & \hspace{-0.22cm}  & n \, + & \hspace{-0.27cm} 2 \hspace{-0.32cm}  & m \, - & \hspace{-0.37cm}  &  c
  & \hspace{0.4cm} & \hspace{0.4cm}    & \alpha'(G) \ge  0.11111 \cdot n + 0.22222 \cdot m - 0.11111 \cdot c \\ \hline
4 &  20 \, \alpha'(G)  & \ge & \hspace{-0.22cm}  & n \, + & \hspace{-0.27cm} 4 \hspace{-0.32cm}  & m \, - & \hspace{-0.37cm}  &  c
  & \hspace{0.4cm} & \hspace{0.4cm}    & \alpha'(G) \ge  0.05000 \cdot n + 0.20000 \cdot m - 0.05000 \cdot c \\ \hline
5 &  55 \, \alpha'(G)  & \ge & \hspace{-0.22cm} 2 \hspace{-0.32cm}  & n \, + & \hspace{-0.27cm} 9 \hspace{-0.32cm}  & m \, - & \hspace{-0.37cm} 2 \hspace{-0.32cm}  &  c
  & \hspace{0.4cm} & \hspace{0.4cm}    & \alpha'(G) \ge  0.03636 \cdot n + 0.16364 \cdot m - 0.03636 \cdot c \\ \hline
6 &  42 \, \alpha'(G)  & \ge & \hspace{-0.22cm}  & n \, + & \hspace{-0.27cm} 6 \hspace{-0.32cm}  & m \, - & \hspace{-0.37cm}  &  c
  & \hspace{0.4cm} & \hspace{0.4cm}    & \alpha'(G) \ge  0.02381 \cdot n + 0.14286 \cdot m - 0.02381 \cdot c \\ \hline
7 &  161 \, \alpha'(G)  & \ge & \hspace{-0.22cm} 3 \hspace{-0.32cm}  & n \, + & \hspace{-0.27cm} 20 \hspace{-0.32cm}  & m \, - & \hspace{-0.37cm} 3 \hspace{-0.32cm}  &  c
  & \hspace{0.4cm} & \hspace{0.4cm}    & \alpha'(G) \ge  0.01863 \cdot n + 0.12422 \cdot m - 0.01863 \cdot c \\ \hline
8 &  72 \, \alpha'(G)  & \ge & \hspace{-0.22cm}  & n \, + & \hspace{-0.27cm} 8 \hspace{-0.32cm}  & m \, - & \hspace{-0.37cm}  &  c
  & \hspace{0.4cm} & \hspace{0.4cm}    & \alpha'(G) \ge  0.01389 \cdot n + 0.11111 \cdot m - 0.01389 \cdot c \\ \hline
9 &  351 \, \alpha'(G)  & \ge & \hspace{-0.22cm} 4 \hspace{-0.32cm}  & n \, + & \hspace{-0.27cm} 35 \hspace{-0.32cm}  & m \, - & \hspace{-0.37cm} 4 \hspace{-0.32cm}  &  c
  & \hspace{0.4cm} & \hspace{0.4cm}    & \alpha'(G) \ge  0.01140 \cdot n + 0.09972 \cdot m - 0.01140 \cdot c \\ \hline
10 &  110 \, \alpha'(G)  & \ge & \hspace{-0.22cm}  & n \, + & \hspace{-0.27cm} 10 \hspace{-0.32cm}  & m \, - & \hspace{-0.37cm}  &  c
  & \hspace{0.4cm} & \hspace{0.4cm}    & \alpha'(G) \ge  0.00909 \cdot n + 0.09091 \cdot m - 0.00909 \cdot c \\ \hline
11 &  649 \, \alpha'(G)  & \ge & \hspace{-0.22cm} 5 \hspace{-0.32cm}  & n \, + & \hspace{-0.27cm} 54 \hspace{-0.32cm}  & m \, - & \hspace{-0.37cm} 5 \hspace{-0.32cm}  &  c
  & \hspace{0.4cm} & \hspace{0.4cm}    & \alpha'(G) \ge  0.00770 \cdot n + 0.08320 \cdot m - 0.00770 \cdot c \\ \hline
\end{array}
\]
\begin{center}
\textbf{Table~2.} Tight lower bounds on the matching number of a graph with maximum degree~$k$ and with no $k$-regular component.
\end{center}

}

The following result presents another lower bound on the matching number when $k \ge 2$ is even.  A proof of Theorem~\ref{match_thm2} is given in Section~\ref{S:match_thm2}.

\begin{thm}\label{match_thm2}
Let $k \ge 2$ be an even integer and let $G$ be any graph of order $n$ and size $m$ and maximum degree $\Delta(G) \le k$.
If no component of $G$ is $k$-regular, then
\[
 \alpha'(G) \ge  \left( \frac{k+2}{k^2+k+2} \right) m  \, - \,  \left( \frac{k-2}{k^2+k+2} \right) n.
\]
\end{thm}

If $G$ is a connected graph of order $n$ with maximum degree at most~$2$, then $G$ is a either a path or a cycle, implying that $\alpha'(G) = \lceil \frac{n-1}{2} \rceil$. Hence, it is only of interest to focus on connected graphs with maximum degree at most~$k$, where $k \ge 3$.  As a consequence of Theorem~\ref{KregEVEN} and Theorem~\ref{match_thm}, we have the following result when $k \ge 4$ is even. A proof of Corollary~\ref{corKeven} is given in Section~\ref{S:corKeven}.

\begin{cor}\label{corKeven}
If $k \ge 4$ is an even integer and $G$ is a connected graph of order $n$, size $m$ and maximum degree $\Delta(G) \le k$, then
\[
\alpha'(G) \ge \frac{n}{k(k+1)} \, + \, \frac{m}{k+1} - \frac{1}{k(k+1)},
\]
unless the following holds.
\\[-26pt]
\begin{enumerate}
\item $G$ is $k$-regular and $n = k+1$, in which case $\alpha'(G) = \frac{n-1}{2} = \frac{n}{k(k+1)} \, + \, \frac{m}{k+1} - \frac{1}{k}$.
\item $G$ is $k$-regular and $n = k+3$, in which case $\alpha'(G) = \frac{n-1}{2} = \frac{n}{k(k+1)} \, + \, \frac{m}{k+1} - \frac{3}{k(k+1)}$.
\end{enumerate}
\end{cor}

As a consequence of Theorem~\ref{KregODD} and Theorem~\ref{match_thm}, we have the following result when $k \ge 3$ is odd. A proof of Corollary~\ref{corKodd} is given in Section~\ref{S:corKodd}.

\begin{cor}\label{corKodd}
If $k \ge 3$ is an odd integer and $G$ is a connected graph of order $n$, size $m$, and with maximum degree $\Delta(G) \le k$, then
\[
\alpha'(G) \ge \left( \frac{k-1}{k(k^2 - 3)} \right) n \, + \, \left( \frac{k^2 - k - 2}{k(k^2 - 3)} \right) m \, - \, \frac{k-1}{k(k^2 - 3)}.
\]
\end{cor}

As a consequence of Theorem~\ref{KregEVEN} and Theorem~\ref{match_thm2}, we have the following result when $k \ge 4$ is even. A proof of Corollary~\ref{corKeven2} is given in Section~\ref{S:corKeven2}.

\begin{cor} \label{corKeven2}
If $k \ge 4$ is an even integer and $G$ is a graph of order $n$, size $m$ and maximum degree $\Delta(G) \le k$, then
\[
 \alpha'(G) \ge  \left( \frac{k+2}{k^2+k+2} \right) m  \, - \,  \left( \frac{k-2}{k^2+k+2} \right) n
\]
unless the following holds.
\\[-26pt]
\begin{enumerate}
\item $G$ is $k$-regular and $n = k+1$, in which case  $\alpha'(G) \ge  \left( \frac{k+2}{k^2+k+2} \right) m  \, - \,  \left( \frac{k-2}{k^2+k+2} \right) n \, - \frac{k+2}{k^2+k+2}$.

\item $G$ is $k$-regular and $n = k+3$, in which case $\alpha'(G) \ge  \left( \frac{k+2}{k^2+k+2} \right) m  \, - \,  \left( \frac{k-2}{k^2+k+2} \right) n  \, - \frac{4}{k^2+k+2}$.

\item $G$ is $4$-regular and $n = 9$, in which case $\alpha'(G) \ge  \left( \frac{k+2}{k^2+k+2} \right) m  \, - \,  \left( \frac{k-2}{k^2+k+2} \right) n  \, - \frac{2}{k^2+k+2}$.
\end{enumerate}
\end{cor}

Theorem~B and Theorem~C follow from Corollaries~\ref{corKeven},~\ref{corKodd}, and~\ref{corKeven2}.

Let $n_i$ denotes the number of vertices of degree~$i$ in a graph $G$. We remark that substituting $k = 3$, $n = n_1 + n_2 + n_3$, and $m = \frac{1}{2}(n_1 + 2n_2 + 3n_3)$ into the lower bound in the statement of Corollary~\ref{corKodd} yields the following result of Haxwell and Scott~\cite{HaSc14}.

\begin{cor}{\rm (\cite{HaSc14})}
If $G$ is a graph with maximum degree $\Delta(G) \le 3$, then
\[
\alpha'(G) \ge  \frac{4}{9} n_3 \, + \, \frac{1}{3} n_2 \, + \, \frac{2}{9} n_1 \, - \, \frac{1}{9}c.
\]
\label{corHaxScott}
\end{cor}

\subsection{Motivation}

Our aim in this paper, is for all $k \ge 3$, to give a complete description of the set $L_k$ of pairs $(\gamma,\beta)$ of real numbers for which there exists a constant $K$  such that
$\alpha'(G) \ge \gamma n + \beta m - K$ holds every connected graph $G$ with maximum degree at most~$k$, where $n$ and $m$ denote the number of vertices and the number of edges, respectively, in $G$. Similar work was done by Chv\'{a}tal and McDiarmid~\cite{ChMc} for the transversal number of a $k$-uniform hypergraph, for $k \ge 2$, in terms of its order and size. In their case, the resulting convex set has infinitely many extreme points. In our case, we show that in contrast to the Chv\'{a}tal-McDiarmid result, our convex set has exactly one extreme point when $k$ is odd and exactly two extreme points when $k$ is even.

Various lower bounds on the matching number for regular graphs have appeared in the literature. For example, Biedl et. al~\cite{Biedl} proved that if $G$ is a cubic graph, then $\alpha'(G) \ge (4n-1)/9$. This result was generalized to regular graphs of  higher degree by Henning and Yeo~\cite{HeYe07} (see also, O and West~\cite{OWest10}). O and West~\cite{OWest11} established lower bounds on the matching number with given edge-connectivity in regular graphs. Cioab\u{a}, Gregory, and Haemers~\cite{CiGrHa09} studied matchings in regular graphs from eigenvalues. Lower bounds on the matching number for general graphs and bipartite graphs were obtained by Jahanbekam and West~\cite{JaWe13}.

Lower bounds on the matching number in subcubic graphs (graphs with maximum degree at most~$3$) were studied by, among others, Henning, L\"{o}wenstein, and  Rautenbach~\cite{HeLoRa12}. Recently, Haxell and Scott~\cite{HaSc14} gave a complete description of the set of triples $(\alpha,\beta,\gamma)$ of real numbers for which there exists a constant $K$ such
that $\alpha'(G) \ge \alpha n_3 + \beta n_2 +
\gamma n_1 - K$ for every connected subcubic graph $G$, where $n_i$ denotes the number of vertices of degree $i$ for each $i \in [3]$. Here, the resulting convex set is shown to be a $3$-dimensional polytope determined by the the intersection of the six half-spaces.

In this paper, we establish a tight lower bound on the matching number of a graph with given maximum degree in terms of its order and size.

For graph theory and terminology, we generally follow~\cite{MHAYbookTD}. In particular, we denote the \emph{degree} of a vertex $v$ in the graph $G$ by $d_G(v)$. The maximum degree among the vertices of $G$ is denoted by $\Delta(G)$. For a subset $S$ of vertices of a graph $G$, we let $G[S]$ denote the subgraph induced by $S$. The number of odd components of a graph $G$ we denote by $\oc(G)$. We use the standard notation $[k] = \{1,2,\ldots,k\}$.

\section{Proof of Theorem~\ref{match_thm}}
\label{S:match_thm1}

Let
\[
\varepsilon_k = \left\{
\begin{array}{lc}
\displaystyle{ \frac{2}{k(k+1)} } & \mbox{for $k \ge 2$ even} \3 \\
\displaystyle{ \frac{2k-2}{k(k^2-3)} } & \mbox{for $k \ge 3$ odd}.
\end{array}
\right.
\]
For $k \ge 2$, let
\[
a_k= \frac{\varepsilon_k}{2} \hspace*{0.5cm} \mbox{and}  \hspace*{0.5cm} b_k = \frac{2 - k \varepsilon_k }{2k}.
\]
We note that for $k \ge 2$ even,
\[
a_k= \frac{1}{k(k+1)} \hspace*{0.5cm} \mbox{and}  \hspace*{0.5cm} b_k = \frac{1}{k+1},
\]
and for $k \ge 3$ odd,
\[
a_k = \frac{k-1}{k(k^2 - 3)} \hspace*{0.5cm} \mbox{and}  \hspace*{0.5cm} b_k = \frac{k^2 - k - 2}{k(k^2 - 3)}.
\]
Further, in both cases,
\[
a_k + b_k = \frac{1}{k},
\]

\noindent
and so $ka_k + kb_k  = 1$. Theorem~\ref{match_thm} can now be restated as follows.

\medskip
\noindent \textbf{Theorem~\ref{match_thm}} \emph{Let $k \ge 4$ be an integer and let $G$ be a graph with $c$ components and of order $n$ and size $m$ and maximum degree $\Delta(G) \le k$. If no component of $G$ is $k$-regular, then $\alpha'(G) \ge a_k n + b_k m - a_k c$. } \\

\noindent \textbf{Proof of Theorem~\ref{match_thm}.} Let $c(G^*)$ denote the number of components of a graph $G^*$.
Define the following five values of a graph $G^*$ and vertex set $X^* \subseteq V(G^*)$.

\begin{description}
 \item[(1):] $\beta_1(G^*,X^*)$ is the number of edges in $G^*[X^*]$.
 \item[(2):] $\beta_2(G^*,X^*)$ is the number of even components in $G^*-X^*$.
 \item[(3):] $\beta_3(G^*,X^*)$ is the number of vertices in $X^*$ with degree less than $k$.
 \item[(4):] $\beta_4(G^*,X^*)$ is the number of components in $G^*-X^*$ that do not have exactly one edge to $X^*$.
 \item[(5):] $\beta_5(G^*,X^*)$ is the number of odd components in $G^*-X^*$ with order between $1$ and $k+1$.
\end{description}

For the sake of contradiction suppose that the theorem is false and that $G$ is a counter example to the theorem.
That is, $G$ has maximum degree at most $k$ and
no component of $G$ is $k$-regular and $\alpha'(G) < a_k \cdot |V(G)| \, + \, b_k \cdot |E(G)| - a_k \cdot c(G)$.
By the Tutte-Berge formula in Theorem~\ref{Berge} we may assume that $G$ and $X$ are chosen such that the following holds.

\begin{equation}
\frac{1}{2} \left( n + |X| - \oc(G - X)  \right) < a_k \cdot |V(G)| \, + \, b_k \cdot |E(G)| - a_k \cdot c(G).
\label{Eq1}
\end{equation}

Furthermore we may assume that $(\beta_1(G,X), \beta_2(G,X), \ldots, \beta_5(G,X))$ is lexicographically minimum of all $G$ and $X$ satisfying the above. We proceed further with the following series of claims.

\begin{unnumbered}{Claim~A}
$\beta_1(G,X)=0$.
\end{unnumbered}
\proof Suppose, to the contrary, that $\beta_1(G,X) \ge 1$, and let $x_1 x_2 \in E(G)$ be arbitrary where $x_1,x_2 \in X$.
Delete the edge $x_1 x_2$ and add a new vertex $u$ and the edges $x_1 u$ and $u x_2$. Let $G'$ be the resulting graph. We note that $c(G')=c(G)$, $|V(G')|=|V(G)|+1$, $|E(G')|=|E(G)|+1$ and $\oc(G'-X) = \oc(G-X)+1$.
Also, no component in $G'$ is $k$-regular (as if any component in $G'$ is $k$-regular, then the corresponding component in $G$
would also be $k$-regular).
The following now holds.

\[
\begin{array}{rcl}
\vspace{0.2cm} \frac{1}{2} \left( |V(G')| + |X| - \oc(G' - X)  \right)  & = &   \frac{1}{2}((|V(G)|+1) + |X| - (\oc(G-X)+1))  \\ \vspace{0.2cm}
& = & \frac{1}{2} \left(|V(G)| +|X| - \oc(G-X)\right) \\ \vspace{0.2cm}
& < & a_k \cdot |V(G)| \, + \, b_k \cdot |E(G)| - a_k \cdot c(G)\\ \vspace{0.2cm}
& < & a_k \cdot |V(G')| \, + \, b_k \cdot |E(G')| - a_k \cdot c(G').
\end{array}
\]

As $\beta_1(G',X) < \beta_1(G,X)$ this contradicts the lexicographical minimality of  $(\beta_1(G,X)$, $\beta_2(G,X)$, $\ldots, \beta_5(G,X))$.~\smallqed

\begin{unnumbered}{Claim~B}
$\beta_2(G,X)=0$.
\end{unnumbered}
\proof Suppose, to the contrary, that $\beta_2(G,X) \ge 1$. Let $C$ be an even component $C$ in $G-X$, and let $u$ be a leaf in some spanning tree of $C$. In this case $C-\{u\}$ is connected. Let $G' = G - u$ and
note that $|V(G')|=|V(G)|-1$, $|E(G')| \ge |E(G)|-k$, $c(G') \le c(G) + k - 1$ (as we delete at most $k-1$ edges not in $C$ incident to $u$)
and $\oc(G'-X)=\oc(G-X)+1$.
The following now holds (as we have shown that $ka_k+kb_k=1$).

\[
\begin{array}{rcl}
\multicolumn{3}{l}{\vspace{0.2cm} \frac{1}{2} \left( |V(G')| + |X| - \oc(G' - X)  \right)  =        \frac{1}{2}((|V(G)|-1) + |X| - (\oc(G-X)+1)) } \\ \vspace{0.2cm}
& = & \frac{1}{2} \left(|V(G)| +|X| - \oc(G-X) \right) -1 \\ \vspace{0.2cm}
& < & a_k \cdot |V(G)| \, + \, b_k \cdot |E(G)| - a_k \cdot c(G) - 1\\ \vspace{0.2cm}
& \le & a_k \cdot (|V(G')|+1) \, + \, b_k \cdot (|E(G')|+k) -  a_k \cdot (c(G')- k + 1) - 1\\ \vspace{0.2cm}
& = & a_k \cdot |V(G')| \, + \, b_k \cdot |E(G')| - a_k \cdot c(G') +(ka_k + k b_k - 1) \\ \vspace{0.2cm}
& = & a_k \cdot |V(G')| \, + \, b_k \cdot |E(G')| - a_k \cdot c(G').
\end{array}
\]

Any component in $G'$ is either a component in $G$ or contains a vertex of degree at most $k-1$ (adjacent to $u$ in $G$), which implies that
no component of $G'$ is $k$-regular. Furthermore as $\beta_1(G',X)=\beta_1(G,X)=0$ and $\beta_2(G',X)<\beta_2(G,X)$ we obtain a
contradiction to the lexicographical minimality of  $(\beta_1(G,X), \beta_2(G,X), \ldots, \beta_5(G,X))$.~\smallqed

\begin{unnumbered}{Claim~C}
$\beta_3(G,X)=0$.
\end{unnumbered}
\proof Suppose, to the contrary, that $\beta_3(G,X) \ge 1$. Let $x$ be a vertex in $X$ with $d_{G}(x)<k$.
Let $G'$ be obtained by adding $s = k-d_{G}(x)$ new vertices $u_1,u_2,\ldots ,u_s$ and the edges $u_1x,u_2x,\ldots,u_sx$ to $G$.
Note that $|V(G')|=|V(G)|+s$ and $|E(G')| = |E(G)|+s$ and $c(G') = c(G)$ and $\oc(G'-X)=\oc(G-X)+s$.
The following now holds.

\[
\begin{array}{rcl}
\multicolumn{3}{l}{\vspace{0.2cm} \frac{1}{2} \left( |V(G')| + |X| - \oc(G' - X)  \right)  =       \frac{1}{2}((|V(G)|+s) + |X| - (\oc(G-X)+s)) } \\ \vspace{0.2cm}
& = & \frac{1}{2} \left(|V(G)| +|X| - \oc(G-X) \right)  \\ \vspace{0.2cm}
& < & a_k \cdot |V(G)| \, + \, b_k \cdot |E(G)| - a_k \cdot c(G) \\ \vspace{0.2cm}
& \le & a_k \cdot (|V(G')|-s) \, + \, b_k \cdot (|E(G')|-s) - a_k \cdot c(G') \\ \vspace{0.2cm}
& < & a_k \cdot |V(G')| \, + \, b_k \cdot |E(G')| - a_k \cdot c(G'). \\
\end{array}
\]

Any component in $G'$ is either a component in $G$ or contains vertices of degree $1< k-1$ (adjacent to $x$ in $G$), which implies that
no component of $G''$ is $k$-regular.
As $\beta_i(G',X) = \beta_i(G,X)$ for $i \in [2]$ and $\beta_3(G',X) < \beta_3(G,X)$  we obtain a contradiction to the lexicographical minimality of
$(\beta_1(G,X), \beta_2(G,X), \ldots, \beta_5(G,X))$.~\smallqed

\begin{unnumbered}{Claim~D}
$\beta_4(G,X)=0$.
\end{unnumbered}
\proof Suppose, to the contrary, that $\beta_4(G,X) \ge 1$. Let $C$ be an odd component in $G-X$ with $s$ edges to $X$, where $s \ne 1$. Let $q = |V(C)|$.

Suppose that $s=0$ and let $G' = G - C$. Suppose further that $q \le k$. In this case, we note that $|E(C)| \le {q \choose 2}$, which implies the following (as $ka_k+kb_k=1$).

\[
\begin{array}{rcl}
\multicolumn{3}{l}{\vspace{0.2cm} \frac{1}{2} \left( |V(G')| + |X| - \oc(G' - X)  \right)  \le        \frac{1}{2}((|V(G)|-q) + |X| - (\oc(G-X)-1)) } \\ \vspace{0.2cm}
& = & \frac{1}{2} \left(|V(G)| +|X| - \oc(G-X) \right) + \frac{1-q}{2} \\ \vspace{0.2cm}
& < & a_k \cdot |V(G)| \, + \, b_k \cdot |E(G)| - a_k \cdot c(G) + \frac{1-q}{2}  \\ \vspace{0.2cm}
& \le & a_k \cdot (|V(G')|+q) \, + \, b_k \cdot \left( |E(G')| + \frac{q(q-1)}{2} \right) - a_k \cdot (c(G')+1) + \frac{1-q}{2}  \\ \vspace{0.2cm}
& = & a_k \cdot |V(G')| \, + \, b_k \cdot |E(G')| - a_k \cdot c(G') - \frac{q-1}{2} \left( 1 - 2a_k - b_k q \right) \\ \vspace{0.2cm}
& \le & a_k \cdot |V(G')| \, + \, b_k \cdot |E(G')| - a_k \cdot c(G') - \frac{q-1}{2} \left( 1 - k \, a_k  - k \, b_k \right) \\ \vspace{0.2cm}
& = & a_k \cdot |V(G')| \, + \, b_k \cdot |E(G')| - a_k \cdot c(G').
\end{array}
\]

As $\beta_i(G',X) = \beta_i(G,X)$ for $i \in [3]$ and $\beta_4(G',X) < \beta_4(G,X)$  we obtain a contradiction to the lexicographical minimality of
$(\beta_1(G,X), \beta_2(G,X), \ldots, \beta_5(G,X))$. Therefore, $q \ge k+1$.

Suppose that $k$ is even. Recall that $a_k = 1/(k(k+1))$ and $b_k=1/(k+1)$. In this case, as $\Delta(G) \le k$ and $C$ is not $k$-regular we have $|E(C)| \le \frac{1}{2}kq - 1$. This implies the following, as $q \ge k+1$.

\[
\begin{array}{rcl}
\multicolumn{3}{l}{\vspace{0.2cm} \frac{1}{2} \left( |V(G')| + |X| - \oc(G' - X)  \right)  =        \frac{1}{2}((|V(G)|-q) + |X| - (\oc(G-X)-1)) } \\ \vspace{0.2cm}
& = & \frac{1}{2} \left(|V(G)| +|X| - \oc(G-X) \right) + \frac{1-q}{2} \\ \vspace{0.2cm}
& < & a_k \cdot |V(G)| \, + \, b_k \cdot |E(G)| - a_k \cdot c(G) + \frac{1-q}{2}  \\ \vspace{0.2cm}
& \le & a_k \cdot (|V(G')|+q) \, + \, b_k \cdot \left( |E(G')| + \frac{kq}{2} -1 \right) - a_k \cdot (c(G')+1) + \frac{1-q}{2}  \\ \vspace{0.2cm}
& = & a_k \cdot |V(G')| \, + \, b_k \cdot |E(G')| - a_k \cdot c(G') - \left( \frac{q-1}{2} - a_k(q-1) - b_k \frac{kq}{2}  + b_k \right) \\ \vspace{0.2cm}
& = & a_k \cdot |V(G')| \, + \, b_k \cdot |E(G')| - a_k \cdot c(G') - \left( \frac{q-1}{2} - \frac{q-1}{k(k+1)} - \frac{kq}{2(k+1)} + \frac{1}{k+1} \right) \\ \vspace{0.2cm}
& = & a_k \cdot |V(G')| \, + \, b_k \cdot |E(G')| - a_k \cdot c(G') - \frac{q}{2(k+1)} \left( (k+1) - \frac{2}{k} - k \right) + \frac{1}{2} - \frac{1}{k(k+1)} - \frac{1}{k+1} \\ \vspace{0.2cm}
& = & a_k \cdot |V(G')| \, + \, b_k \cdot |E(G')| - a_k \cdot c(G') - \frac{q}{2(k+1)} \left( 1 - \frac{2}{k} \right) + \frac{k(k+1) - 2 - 2k}{2k(k+1)}  \\ \vspace{0.2cm}
& \le & a_k \cdot |V(G')| \, + \, b_k \cdot |E(G')| - a_k \cdot c(G') - \frac{1}{2} \left( \frac{k-2}{k} \right) + \frac{k^2-k-2}{2k(k+1)}  \\ \vspace{0.2cm}
& = & a_k \cdot |V(G')| \, + \, b_k \cdot |E(G')| - a_k \cdot c(G') + \frac{-(k-2)(k+1)+k^2-k-2}{2k(k+1)} \\ \vspace{0.2cm}
& = & a_k \cdot |V(G')| \, + \, b_k \cdot |E(G')| - a_k \cdot c(G').
\end{array}
\]

As $\beta_i(G',X) = \beta_i(G,X)$ for $i \in [3]$ and $\beta_4(G',X) < \beta_4(G,X)$  we obtain a contradiction to the lexicographical minimality of
$(\beta_1(G,X), \beta_2(G,X), \ldots, \beta_5(G,X))$.
Therefore, $k$ is odd, and so $a_k = (k-1)/(k(k^2 - 3))$ and $b_k=(k^2 - k - 2)/(k(k^2 - 3))$. In this case, as $\Delta(G) \le k$ and $C$ is not $k$-regular we have $|E(C)| \le \frac{1}{2}(k q - 1)$. Further, since $C$ is an odd component, we note that $q$ is odd and $q \ge k+2$. This implies the following.

\[
\begin{array}{rcl}
\multicolumn{3}{l}{\vspace{0.2cm} \frac{1}{2} \left( |V(G')| + |X| - \oc(G' - X)  \right)  =        \frac{1}{2}((|V(G)|-q) + |X| - (\oc(G-X)-1)) } \2 \\ 
& = & \frac{1}{2} \left(|V(G)| +|X| - \oc(G-X) \right) + \frac{1-q}{2} \2 \\ 
& < & a_k \cdot |V(G)| \, + \, b_k \cdot |E(G)| - a_k \cdot c(G) + \frac{1-q}{2}  \2 \\ 
& \le & a_k \cdot (|V(G')|+q) \, + \, b_k \cdot \left( |E(G')| + \frac{1}{2}(kq - 1) \right) - a_k \cdot (c(G')+1) + \frac{1-q}{2}  \2 \\ 
& = & a_k \cdot |V(G')| \, + \, b_k \cdot |E(G')| - a_k \cdot c(G') - \left( \frac{q-1}{2} - (q-1)a_k - \frac{1}{2}(kq - 1)b_k  \right) \2 \\ 
& = & a_k \cdot |V(G')| \, + \, b_k \cdot |E(G')| - a_k \cdot c(G') - \left( \frac{q-1}{2} - \frac{(k-1)(q-1)}{k(k^2 - 3)} - \frac{(kq - 1)(k^2 - k - 2)}{2k(k^2 - 3)}  \right) \2 \\ 
& = & a_k \cdot |V(G')| \, + \, b_k \cdot |E(G')| - a_k \cdot c(G') - \1 \\
& & \hspace*{0.5cm} \frac{q}{2k(k^2 - 3)} \left( k(k^2 - 3) - 2(k-1) - k(k^2 - k - 2) \right) + \1 \\
& & \hspace*{1cm} \frac{1}{2k(k^2 - 3)} \left( k(k^2 - 3) - 2(k-1) - (k^2 - k - 2) \right) \2 \\
& = & a_k \cdot |V(G')| \, + \, b_k \cdot |E(G')| - a_k \cdot c(G') - \frac{q(k^2-3k+2)}{2k(k^2 - 3)} + \frac{k^3 - k^2 - 4k + 4}{2k(k^2 - 3)}  \2
\\
& \le & a_k \cdot |V(G')| \, + \, b_k \cdot |E(G')| - a_k \cdot c(G') - \frac{(k+2)(k^2-3k+2)}{2k(k^2 - 3)} + \frac{k^3 - k^2 - 4k + 4}{2k(k^2 - 3)}  \2
\\ 
& = & a_k \cdot |V(G')| \, + \, b_k \cdot |E(G')| - a_k \cdot c(G').
\end{array}
\]

As $\beta_i(G',X) = \beta_i(G,X)$ for $i \in [3]$ and $\beta_4(G',X) < \beta_4(G,X)$  we obtain a contradiction to the lexicographical minimality of
$(\beta_1(G,X), \beta_2(G,X), \ldots, \beta_5(G,X))$.
Therefore, $s \ge 2$.

Let $x_1u_1, x_2u_2, \ldots , x_s u_s$ be distinct edges from $x_i \in X$ to $u_i \in V(C)$ for $i \in [s]$.
Let $G'$ be obtained from $G$ by adding $s-1$ new vertices $w_2,w_3,\ldots,w_s$, deleting the edges $x_2u_2,x_3u_3,\ldots,x_su_s$ and adding the edges
$x_2w_2,x_3w_3,\ldots,x_sw_s$.  It is not difficult to see that no component of $G'$ is $k$-regular.
Note that $|V(G')|=|V(G)|+s-1$ and $|E(G')| = |E(G)|$ and $c(G') \le c(G)+s-1$ and $\oc(G'-X) = \oc(G-X)+s-1$.
The following now holds.

\[
\begin{array}{rcl}
\multicolumn{3}{l}{\vspace{0.2cm} \frac{1}{2} \left( |V(G')| + |X| - \oc(G' - X)  \right)  =        \frac{1}{2}((|V(G)|+s-1) + |X| - (\oc(G-X)+s-1)) } \\ \vspace{0.2cm}
& = & \frac{1}{2} \left(|V(G)| +|X| - \oc(G-X) \right)  \\ \vspace{0.2cm}
& < & a_k \cdot |V(G)| \, + \, b_k \cdot |E(G)| - a_k \cdot c(G) \\ \vspace{0.2cm}
& \le & a_k \cdot (|V(G')|-s+1) \, + \, b_k \cdot |E(G')| - a_k \cdot (c(G')-s+1) \\ \vspace{0.2cm}
& = & a_k \cdot |V(G')| \, + \, b_k \cdot |E(G')| - a_k \cdot c(G').
\end{array}
\]

As $\beta_i(G',X) = \beta_i(G,X)$ for $i \in [3]$ and $\beta_4(G',X) < \beta_4(G,X)$  we obtain a contradiction to the lexicographical minimality of
$(\beta_1(G,X), \beta_2(G,X), \ldots, \beta_5(G,X))$.~\smallqed

\begin{unnumbered}{Claim~E}
$\beta_5(G,X)=0$.
\end{unnumbered}
\proof Suppose, to the contrary, that $\beta_5(G,X) \ge 1$. Let $C$ be an odd component $C$ in $G-X$ with $1 < |V(C)| < k+1$. As $\beta_4(G,X)=0$ there is exactly one edge from $C$ to $X$. Let $xc$ be the edge with $x \in X$ and $c \in C$.
Let $r = |V(C)|-1$.

Let $G' = G - (V(C)\setminus\{c\})$. That is, $G'$ is obtained from $G$ by removing all vertices of $C$, except $c$, from $G$.
Note that $|V(G')|=|V(G)|-r$ and $|E(G')| = |E(G)| - |E(C)|$ and $c(G') = c(G)$ and $\oc(G'-X) = \oc(G-X)$.
It is not difficult to see that no component of $G'$ is $k$-regular.
The following now holds.

\[
\begin{array}{rcl}
\multicolumn{3}{l}{\vspace{0.2cm} \frac{1}{2} \left( |V(G')| + |X| - \oc(G' - X)  \right)  =   \frac{1}{2}((|V(G)|-r) + |X| - (\oc(G-X))) } \\ \vspace{0.2cm}
& = & \frac{1}{2} \left(|V(G)| +|X| - \oc(G-X) \right) - \frac{r}{2} \\ \vspace{0.2cm}
& < & a_k \cdot |V(G)| \, + \, b_k \cdot |E(G)| - a_k \cdot c(G) - \frac{r}{2}  \\ \vspace{0.2cm}
& = & a_k \cdot (|V(G')|+r) \, + \, b_k \cdot (|E(G')|+|E(C)| ) - a_k \cdot c(G') - \frac{r}{2} \\ \vspace{0.2cm}
& = & a_k \cdot |V(G')| \, + \, b_k \cdot |E(G')| - a_k \cdot c(G') - \left( \frac{r}{2} - a_k r - b_k |E(C)| \right).
\end{array}
\]
We will now evaluate $\frac{r}{2} - a_k r - b_k |E(C)|$. Note that $|E(C)| \le \frac{r(r+1)}{2}$ as $|V(C)|=r+1$ and $r \le k-1$.
If $k$ is even, then $a_k = 1/(k(k+1))$ and $b_k=1/(k+1)$, which implies the following (as $0<r<k$ and $k \ge 2$).

\[
\begin{array}{rcl}
\vspace{0.2cm} \frac{r}{2} - a_k r - b_k |E(C)| & \ge & \frac{r}{2} - \frac{1}{k(k+1)} r - \frac{1}{k+1} \times \frac{r(r+1)}{2}  \\ \vspace{0.2cm}
& = & \frac{r}{2k(k+1)} \left( k(k+1) - 2 - k(r+1) \right) \\ \vspace{0.2cm}
& = & \frac{r}{2k(k+1)} \left( k(k-r) - 2 \right) \\ \vspace{0.2cm}
& \ge & 0.
\end{array}
\]

If $k$ is odd, then $a_k = (k-1)/(k(k^2 - 3))$ and $b_k=(k^2 - k - 2)/(k(k^2 - 3))$, which implies the following (as $0<r<k$ and $k \ge 3$).

\[
\begin{array}{rcl}
\vspace{0.2cm} \frac{r}{2} - a_k r - b_k |E(C)| & \ge & \frac{r}{2} - \frac{(k-1)r}{k(k^2 - 3)} - \frac{r(r+1)(k^2 - k - 2)}{2k(k^2 - 3)}  \\ \vspace{0.2cm}
& = & \frac{r}{2k(k^2 - 3)} \left( k(k^2 - 3) - 2(k-1) - (k^2 - k - 2)(r+1) \right) \\ \vspace{0.2cm}
& \ge & \frac{r}{2k(k^2 - 3)} \left( k^3 - 5k + 2 - (k^2 - k - 2)k \right) \\ \vspace{0.2cm}
& \ge & \frac{r}{2k(k^2 - 3)} \left( k^2 - 3k + 2 \right) \\ \vspace{0.2cm}
& \ge & \frac{r(k-1)(k-2)}{2k(k^2 - 3)} \\ \vspace{0.2cm}
& > & 0.
\end{array}
\]

In both cases, $\frac{r}{2} - a_k r - b_k |E(C)| \ge 0$. This implies the following.

\[
\frac{1}{2} \left( |V(G')| + |X| - \oc(G' - X)  \right) < a_k \cdot |V(G')| \, + \, b_k \cdot |E(G')| - a_k \cdot c(G').
\]

As $\beta_i(G',X) = \beta_i(G,X)$ for $i \in [4]$ and $\beta_5(G',X) < \beta_5(G,X)$  we obtain a contradiction to the lexicographical minimality of
$(\beta_1(G,X), \beta_2(G,X), \ldots, \beta_5(G,X))$.~\smallqed


\bigskip
We now return to the proof of Theorem~\ref{match_thm}. Let $\cC$ be the set of all components of $G$, and so $|\cC(G)| = c(G)$.

\begin{unnumbered}{Claim~F}
If $C \in \cC$, then
\[
\frac{1}{2} \left(|V(C)| + |X_C| - \oc(C-X_C) \right)  \ge a_k \cdot |V(C)| \, + \, b_k \cdot |E(C)| - a_k.
\]
\end{unnumbered}
\proof Let $C \in \cC$, and let $X_C = X \cap V(C)$. By Claim~A, $\beta_1(G,X)=0$, and so there is no edge in $G[X_C]$. By Claim~D, $\beta_4(G,X)=0$, and so every component in $C-X_C$ has exactly one edge to $X_C$. Since $C$ is connected, we therefore note that $|X_C|=1$.
Let $X_c=\{x\}$. By Claim~C, $\beta_3(G,X)=0$, and so $d_G(x)=d_C(x)=k$.
Let $n_i$ denote that number of components in $C-X_C$ of order $i$, and so
\[
\oc(C-X_C)=k= \sum_{i=1}^{\infty} n_i.
\]

By Claim~B, $\beta_2(G,X) = 0$ and by Claim~E, $\beta_5(G,X)=0$, implying that all $n_i$ are zero except possibly if $i=1$ or $i \ge k+1$ and $i$ is odd. We note that the expression $\frac{1}{2} (|V(C)| + |X_C| - \oc(C-X_C))$ can be written in terms of $k$ and $\varepsilon_k$ as follows.

\[
\begin{array}{rcl}
\multicolumn{3}{l}{\vspace{0.2cm} \displaystyle{ \frac{1}{2} \left(|V(C)| + |X_C| - \oc(C-X_C) \right) }} \\ \vspace{0.2cm}
& = &
\displaystyle{ \frac{1}{2} \left(|X_C| + \left( \sum_{i=1}^{\infty} i \cdot n_i \right) + |X_C| - \oc(G-X_G) \right) } \\ \vspace{0.2cm}
& = & \displaystyle{ \frac{1}{2} \left(2 + \left( \sum_{i=1}^{\infty} i \cdot n_i \right) - (1-\varepsilon_k) \, \oc(G-X_G) - \varepsilon_k \, \oc(G-X_G) \right) } \\ \vspace{0.2cm}
& = & \displaystyle{ \frac{1}{2} \left(2 + \left( \sum_{i=1}^{\infty} i \cdot n_i \right) - (1-\varepsilon_k) \, \left( \sum_{i=1}^{\infty} n_i \right) - \varepsilon_k \, k \right) } \\ \vspace{0.2cm}
& = & \displaystyle{ 1 - \left( \frac{\varepsilon_k}{2} \right) k + \sum_{i=1}^{\infty} \left( i - 1 + \varepsilon_k \right) \frac{n_i}{2} } \\ \vspace{0.2cm}
& = & \displaystyle{ 1 - \left( \frac{\varepsilon_k}{2} \right) k + \left( \frac{\varepsilon_k}{2} \right) n_1 + \sum_{i=k+1}^{\infty} \left( i - 1 + \varepsilon_k \right) \frac{n_i}{2} }. \\ \vspace{0.2cm}
\end{array}
\]
Recall that
\[
a_k = \frac{\varepsilon_k}{2} \hspace*{0.5cm} \mbox{and}  \hspace*{0.5cm} b_k = \frac{2 - k \varepsilon_k }{2k}.
\]

Hence, the expression $\frac{1}{2} (|V(C)| + |X_C| - \oc(C-X_C))$ can be written as follows.
\begin{equation}
\frac{1}{2} \left(|V(C)| + |X_C| - \oc(C-X_C) \right) = kb_k + a_kn_1 + \sum_{i=k+1}^{\infty} \left( i - 1 + \varepsilon_k \right) \frac{n_i}{2}.
\label{Eq2}
\end{equation}
We consider two cases, depending on the parity of $k$.

\emph{Case~1. $k$ is even.} For all components, $C^*$, in $C - X_C$ we note that $C^*$ is not $k$-regular. Further, if $C^*$ has order $r$, then since $k$ is even, $|E(C^*)| \le (rk-2)/2$. Thus,

\[
\begin{array}{rcl}  
\multicolumn{3}{l}{\vspace{0.2cm} \displaystyle{ a_k \cdot |V(C)| + b_k \cdot |E(C)| - a_k  }} \\ \vspace{0.2cm}
& \le &
  \displaystyle{ a_k \left( |X_C|+ \sum_{i=1}^{\infty} i \cdot n_i \right) + b_k \left( k + \sum_{i=k+1}^{\infty} \frac{ki-2}{2} \right)  - a_k } \\  \vspace{0.2cm}
& = & \displaystyle{ a_k \left( 1 + n_1 + \sum_{i=k+1}^{\infty} i \cdot n_i \right) + b_k \left( k + \sum_{i=k+1}^{\infty} \frac{ki-2}{2} \right)  - a_k } \\  \vspace{0.2cm}
& = & \displaystyle{ k b_k + a_k n_1 + \sum_{i=k+1}^{\infty} \left( a_i \cdot i + b_k \left( \frac{ki-2}{2} \right) \right) n_i }.
\end{array}
\]

\vspace{-1.9cm}

\begin{equation}
\label{EqX3}
\end{equation}

\vspace{0.6cm}


By Equation~(\ref{Eq2}) and Inequality~(\ref{EqX3}), we note that the desired result follows if the following is true for all $i \ge k+1$.

\[
 \begin{array}{ccrcl}
 & & \displaystyle{ \frac{1}{2} \left( i - 1 + \varepsilon_k \right) } & \ge & \displaystyle{ a_k \cdot i + b_k \left( \frac{ki-2}{2} \right) } \\
\Updownarrow & & & & \\
 & & \displaystyle{ \frac{1}{2} \left( i - 1 + \varepsilon_k \right) } & \ge & \displaystyle{ \left(\frac{\varepsilon_k}{2}\right) i + \left(\frac{1}{k} - \frac{\varepsilon_k}{2} \right) \left( \frac{ki-2}{2} \right) } \\
\Updownarrow & & & & \\
 & & \displaystyle{ 2k \left( i - 1 + \varepsilon_k \right) } & \ge & \displaystyle{ 2 k \varepsilon_k i + \left( 2 - k \varepsilon_k \right) ( ki-2) } \\
\Updownarrow & & & & \\
  & & \displaystyle{  2ki - 2k + 2k \varepsilon_k } & \ge & \displaystyle{ 2 k \varepsilon_k i +  2ki - i k^2 \varepsilon_k - 4 + 2k \varepsilon_k  } \\
\Updownarrow & & & & \\
  & & \displaystyle{ i k \varepsilon_k (k-2) } & \ge & \displaystyle{ 2 k  -  4  } \\
\Uparrow & & & & \\
  & & \displaystyle{ (k+1) k \varepsilon_k (k-2) } & \ge & \displaystyle{ 2 (k-2)  } \\
\Updownarrow & & & & \\
  & & \displaystyle{  \varepsilon_k } & \ge & \displaystyle{ \frac{2}{k(k+1)}  } \\
\end{array}
\]

The above clearly holds since in this case when $k$ is even, $\varepsilon_k = 2/(k(k+1)$.

\medskip
\emph{Case~2. $k \ge 3$ is odd.} In this case, we note that all $n_i$ are zero except possibly if $i=1$ or $i \ge k+2$ and $i$ is odd.
In particular, we note that in Equation~(\ref{Eq2}) the term $n_{k+1} = 0$.
For all components, $C^*$, in $C - X_C$ we note that $C^*$ is not $k$-regular. Further, if $C^*$ has order $r$, then since $k$ is odd, $|E(C^*)| \le (rk-1)/2$. Thus,

\[
\begin{array}{rcl}
\multicolumn{3}{l}{\vspace{0.2cm} \displaystyle{ a_k \cdot |V(C)| + b_k \cdot |E(C)| - a_k  }} \\ \vspace{0.2cm}
& \le &
  \displaystyle{ a_k \left( |X_C|+ \sum_{i=1}^{\infty} i \cdot n_i \right) + b_k \left( k + \sum_{i=k+2}^{\infty} \left( \frac{ki-1}{2} \right) n_i \right)  - a_k } \\  \vspace{0.2cm}
& = & \displaystyle{ a_k \left( 1 + n_1 + \sum_{i=k+2}^{\infty} i \cdot n_i \right) + b_k \left( k + \sum_{i=k+2}^{\infty} \left( \frac{ki-1}{2} \right) n_i \right)  - a_k } \\  \vspace{0.2cm}
& = & \displaystyle{ k b_k + a_k n_1 + \sum_{i=k+2}^{\infty} \left( a_k \cdot i + b_k \left( \frac{ki-1}{2} \right) \right) n_i }. \\  \vspace{0.2cm}
\end{array}
\]

\vspace{-2.1cm}

\begin{equation}
\label{EqX4}
\end{equation}

\vspace{0.7cm}

By Equation~(\ref{Eq2}) and Inequality~(\ref{EqX4}), we note that the desired result follows if the following is true for all $i \ge k+2$.

\[
 \begin{array}{ccrcl}
 & & \displaystyle{ \frac{1}{2} \left( i - 1 + \varepsilon_k \right) } & \ge & \displaystyle{ a_k \cdot i + b_k \left( \frac{ki-1}{2} \right) } \\
\Updownarrow & & & & \\
 & & \displaystyle{ \frac{1}{2} \left( i - 1 + \varepsilon_k \right) } & \ge & \displaystyle{\frac{\varepsilon_k}{2} i + \left(\frac{1}{k} - \frac{\varepsilon_k}{2} \right) \left( \frac{ki-1}{2} \right) } \\
\Updownarrow & & & & \\
 & & \displaystyle{ 2k \left( i - 1 + \varepsilon_k \right) } & \ge & \displaystyle{ 2 k \varepsilon_k i + \left( 2 - k \varepsilon_k \right) ( ki-1) } \\
\Updownarrow & & & & \\
  & & \displaystyle{  2ki - 2k + 2k \varepsilon_k } & \ge & \displaystyle{ 2 k \varepsilon_k i +  2ki - i k^2 \varepsilon_k - 2 + k \varepsilon_k  } \\
\Updownarrow & & & & \\
  & & \displaystyle{ i k \varepsilon_k (k-2) } & \ge & \displaystyle{ 2 k  -  k \varepsilon_k - 2  } \\
\Uparrow & & & & \\
  & & \displaystyle{ (k+2) k \varepsilon_k (k-2) } & \ge & \displaystyle{ 2 k  -  k \varepsilon_k - 2  } \\
\Updownarrow & & & & \\
  & & \displaystyle{  k \varepsilon_k (k^2-4+1) } & \ge & \displaystyle{ 2 k  - 2  } \\
\Updownarrow & & & & \\
  & & \displaystyle{  \varepsilon_k } & \ge & \displaystyle{ \frac{2k-2}{k(k^2-3)}  } \\
\end{array}
\]

The above clearly holds  in this case when $k$ is odd, $\varepsilon_k = (2k-2)/(k(k^2-3))$. This completes the proof of Claim~F.~\smallqed

Applying Claim~F to each component $C$ in $\cC$, the following holds.
\[
\begin{array}{lcl}
\displaystyle{ \frac{1}{2} \left( n + |X| - \oc(G - X)  \right) } & = & \displaystyle{ \sum_{C \in \cC} \frac{1}{2} \left( |V(C)| + |X_C| - \oc(C-X_C) \right) } \2 \\
& \ge & \displaystyle{ \sum_{C \in \cC} (a_k \cdot |V(C)| \, + \, b_k \cdot |E(C)| - a_k) } \2 \\
& = & a_k \cdot |V(G)| \, + \, b_k \cdot |E(G)| - a_k \cdot c(G),
\end{array}
\]
which contradicts Inequality~(\ref{Eq1}), thereby proving the theorem.~\qed

\section{Proof of Theorem~\ref{match_thm2}}
\label{S:match_thm2}

For all even $k \ge 4$, let
\[
a_k= \frac{k-2}{k^2+k+2} \hspace*{0.5cm} \mbox{and}  \hspace*{0.5cm} b_k = \frac{k+2}{k^2+k+2}.
\]

We note that for $k \ge 2$ even,
\[
k \, b_k \, - \, a_k \, = \frac{k(k+2)-(k-2)}{k^2+k+2} \, = \, 1.
\]

\noindent
Theorem~\ref{match_thm2} can now be restated as follows.

\medskip
\noindent \textbf{Theorem~\ref{match_thm2}} \emph{ Let $k \ge 2$ be an even integer and let $G$ be any graph of order $n$ and size $m$ and maximum degree $\Delta(G) \le k$.
If no component of $G$ is $k$-regular, then  $\alpha'(G) \ge  b_k \cdot m  - a_k \cdot n$.}

\noindent \textbf{Proof of Theorem~\ref{match_thm2}.} Let $K_{k+1}-e$ denote the complete graph on $k+1$ vertices after removing one edge.
Define the following five values of a graph $G^*$ and vertex set $X^* \subseteq V(G^*)$.

\begin{description}
 \item[(1):] $\xi_1(G^*,X^*)$ is the number of edges in $G^*[X^*]$.
 \item[(2):] $\xi_2(G^*,X^*)$ is the number of even components in $G^*-X^*$.
 \item[(3):] $\xi_3(G^*,X^*)$ is the number of vertices in $X^*$ with degree less than $k$.
 \item[(4):] $\xi_4(G^*,X^*)$ is the number of components in $G^*-X^*$ that do not have exactly one edge to $X^*$.
 \item[(5):] $\xi_5(G^*,X^*)$ is the number of odd components in $G^*-X^*$ not isomorphic to $K_{k+1}-e$.
\end{description}

For the sake of contradiction suppose that the theorem is false and that $G$ is a counter example to the theorem.
That is, $G$ has maximum degree at most $k$ and
no component of $G$ is $k$-regular and $\alpha'(G) < b_k \cdot |E(G)| - a_k \cdot |V(G)|$.
By Theorem~\ref{Berge} (the Tutte-Berge formula) we may assume that $G$ and $X$ are chosen such that the following holds.

\begin{equation}
\frac{1}{2} \left( n + |X| - \oc(G - X)  \right) < b_k \cdot |E(G)| \, - \, a_k \cdot |V(G)|.
\label{Eq1b}
\end{equation}

Furthermore we may assume that $(\xi_1(G,X), \xi_2(G,X), \ldots, \xi_5(G,X))$ is lexicographically minimum of all $G$ and $X$ satisfying the above.
Note that if $G$ is not connected, then one of the components of $G$ is also a counter example to the theorem and this component either has the same
value of $(\xi_1(G,X), \xi_2(G,X), \ldots, \xi_5(G,X))$ or smaller. We may therefore assume that $G$ is connected, for otherwise we consider the before mentioned component of $G$. We proceed further with the following series of claims.

\newpage
\begin{unnumbered}{Claim~I}
$\xi_1(G,X)=0$.
\end{unnumbered}
\proof Suppose, to the contrary, that $\xi_1(G,X) \ge 1$, and let $x_1 x_2 \in E(G)$ be arbitrary where $x_1,x_2 \in X$.
Delete the edge $x_1 x_2$ and add a new vertex $u$ and the edges $x_1 u$ and $u x_2$. Let $G'$ be the resulting graph.
We note that $|V(G')|=|V(G)|+1$, $|E(G')|=|E(G)|+1$ and $\oc(G'-X) = \oc(G-X)+1$.
Also, no component in $G'$ is $k$-regular (as if any component in $G'$ is $k$-regular, then the corresponding component in $G$ would also be $k$-regular).
The following now holds, as $b_k>a_k$.

\[
\begin{array}{rcl}
\vspace{0.2cm} \frac{1}{2} \left( |V(G')| + |X| - \oc(G' - X)  \right)  & = &   \frac{1}{2}((|V(G)|+1) + |X| - (\oc(G-X)+1))  \\ \vspace{0.2cm}
& = & \frac{1}{2} \left(|V(G)| +|X| - \oc(G-X)\right) \\ \vspace{0.2cm}
& < & b_k \cdot |E(G)| - a_k \cdot |V(G)| \\ \vspace{0.2cm}
& < &  b_k (|E(G)|+1) - a_k (|V(G)|+1) \\ \vspace{0.2cm}
& = &  b_k \cdot |E(G')| - a_k \cdot |V(G')|.
\end{array}
\]

As $\xi_1(G',X) < \xi_1(G,X)$ this contradicts the lexicographical minimality of  $(\xi_1(G,X)$, $\xi_2(G,X)$, $\ldots, \xi_5(G,X))$.~\smallqed

\begin{unnumbered}{Claim~II}
$\xi_2(G,X)=0$.
\end{unnumbered}
\proof Suppose, to the contrary, that $\xi_2(G,X) \ge 1$.
Let $C$ be an even component $C$ in $G-X$, and let $u$ be a leaf in some spanning tree of $C$. In this case $C-\{u\}$ is connected. Let $G' = G - u$ and
note that $|V(G')|=|V(G)|-1$, $|E(G')| \ge |E(G)|-k$  and $\oc(G'-X)=\oc(G-X)+1$.
The following now holds (as we have shown that $kb_k-a_k=1$).

\[
\begin{array}{rcl}
\multicolumn{3}{l}{\vspace{0.2cm} \frac{1}{2} \left( |V(G')| + |X| - \oc(G' - X)  \right)  =        \frac{1}{2}((|V(G)|-1) + |X| - (\oc(G-X)+1)) } \\ \vspace{0.2cm}
& = & \frac{1}{2} \left(|V(G)| +|X| - \oc(G-X) \right) -1 \\ \vspace{0.2cm}
& < & b_k \cdot |E(G)| \, - \, a_k \cdot |V(G)|  - 1\\ \vspace{0.2cm}
& \le & b_k \cdot (|E(G')|+k)  \, - \, a_k \cdot (|V(G')|+1) - 1\\ \vspace{0.2cm}
& = &  b_k \cdot |E(G')| \, - \, a_k \cdot |V(G')|.
\end{array}
\]

Any component in $G'$ is either a component in $G$ or contains a vertex of degree at most $k-1$ (adjacent to $u$ in $G$), which implies that
no component of $G'$ is $k$-regular. Furthermore as $\xi_1(G',X)=\xi_1(G,X)=0$ and $\xi_2(G',X)<\xi_2(G,X)$ we obtain a
contradiction to the lexicographical minimality of  $(\xi_1(G,X), \xi_2(G,X), \ldots, \xi_5(G,X))$.~\smallqed

\newpage
\begin{unnumbered}{Claim~III}
$\xi_3(G,X)=0$.
\end{unnumbered}
\proof Suppose, to the contrary, that $\xi_3(G,X) \ge 1$. Let $x$ be a vertex in $X$ with $d_{G}(x)<k$.
Let $G'$ be obtained by adding $s = k-d_{G}(x)$ new vertices $u_1,u_2,\ldots ,u_s$ and the edges $u_1x,u_2x,\ldots,u_sx$ to $G$.
Note that $|V(G')|=|V(G)|+s$ and $|E(G')| = |E(G)|+s$ and $\oc(G'-X)=\oc(G-X)+s$.
The following now holds, as $b_k > a_k$.


\[
\begin{array}{rcl}  \vspace{0.2cm}
\frac{1}{2} \left( |V(G')| + |X| - \oc(G' - X)  \right)         & = &      \frac{1}{2}((|V(G)|+s) + |X| - (\oc(G-X)+s))  \\ \vspace{0.2cm}
& = & \frac{1}{2} \left(|V(G)| +|X| - \oc(G-X) \right)  \\ \vspace{0.2cm}
& < &  b_k \cdot |E(G)| \, - \, a_k \cdot |V(G)| \\ \vspace{0.2cm}
& = & b_k \cdot (|E(G')|-s) \, - \, a_k \cdot (|V(G')|-s)  \\ \vspace{0.2cm}
& = & b_k \cdot |E(G')| \, - \, a_k \cdot |V(G')| \, - \, s(b_k-a_k) \\ \vspace{0.2cm}
& < & b_k \cdot |E(G')| \, - \, a_k \cdot |V(G')|.
\end{array}
\]

Any component in $G'$ is either a component in $G$ or contains vertices of degree $1< k$ (adjacent to $x$ in $G$), which implies that
no component of $G'$ is $k$-regular.
As $\xi_i(G',X) = \xi_i(G,X)$ for $i \in [2]$ and $\xi_3(G',X) < \xi_3(G,X)$  we obtain a contradiction to the lexicographical minimality of
$(\xi_1(G,X), \xi_2(G,X), \ldots, \xi_5(G,X))$.~\smallqed


\begin{unnumbered}{Claim~IV}
$\xi_4(G,X)=0$.
\end{unnumbered}
\proof Suppose, to the contrary, that $\xi_4(G,X) \ge 1$. Let $C$ be an odd component in $G-X$ with $s$ edges to $X$, where $s \ne 1$. Let $q = |V(C)|$.

Suppose that $s=0$, which as $G$ is connected (which was proved before Claim~I) implies that $|X|=0$. By Claim~II we note that $q$ is odd.
By Equation~(\ref{Eq1b}) the following holds.

\begin{equation} \label{EqN}
\frac{q-1}{2} = \frac{1}{2} \left( n + |X| - \oc(G - X)  \right) < b_k \cdot |E(G)| \, - \, a_k \cdot |V(G)|.
\end{equation}

Suppose that $q \le k$. In this case, we note that $|E(C)| \le {q \choose 2}$, which implies the following.

\[
        \frac{q-1}{2}
<    b_k \cdot |E(G)| \, - \, a_k \cdot |V(G)|
\le    b_k \cdot \frac{q(q-1)}{2} \, - \, a_k \cdot  q
 =      q \left( \left( \frac{q-1}{2} \right) b_k \, - \, a_k \right).
\]

If $q=1$, then $b_k \cdot |E(G)| \, - \, a_k \cdot |V(G)| \le - a_k < 0 = (q-1)/2$, a contradiction.
By Claim~II we note that $q$ is odd, and so $q \ge 3$, which implies that $b_k \frac{(q-1)}{2} - a_k >0$.
Therefore the following holds.

\[
\begin{array}{rcl}  \vspace{0.2cm}
\frac{q-1}{2} %
& < &  q \left( b_k \frac{(q-1)}{2} \, - \, a_k \right)   \\ \vspace{0.2cm}
& \le &  k \left( b_k \frac{(q-1)}{2} \, - \, a_k \right)  \\ \vspace{0.2cm}
& = &  k \left( \frac{k+2}{k^2+k+2} \times \frac{(q-1)}{2} \, - \, \frac{k-2}{k^2+k+2} \right)  \\ \vspace{0.2cm}
& = &  \frac{q-1}{2} + \left( \frac{q-1}{2} \right) \left( \frac{k(k+2)}{k^2+k+2} -1 \right) \, - \, \frac{k(k-2)}{k^2+k+2}  \\ \vspace{0.2cm}
& = &  \frac{q-1}{2} + \left( \frac{q-1}{2} \right) \left( \frac{k-2}{k^2+k+2} \right) \, - \, k \left( \frac{k-2}{k^2+k+2} \right) \\ \vspace{0.2cm}
& = &  \frac{q-1}{2} + \left( \frac{k-2}{k^2+k+2} \right) \left( \frac{q-1}{2} - k \right)  \\ \vspace{0.2cm}
& < &  \frac{q-1}{2}.
\end{array}
\]

This is clearly a contradiction, which implies that, $q \ge k+1$. Therefore, since $G$ is not $k$-regular, $|E(G)| \le \frac{kq}{2} - 1$ and

\[
\begin{array}{rcl} \vspace{0.2cm}
\frac{q-1}{2}
& < &  b_k \left( \frac{qk-2}{2} \right) \, - \, a_k \cdot q   \\ \vspace{0.2cm}
& = &  \frac{q}{2} \left( kb_k -a_k \right) \, - \, \frac{a_kq}{2} \, - \, b_k  \\ \vspace{0.2cm}
& = &  \frac{q}{2}  \, - \, \frac{a_kq}{2} \, - \, b_k  \\ \vspace{0.2cm}
& \le &  \frac{q}{2}  \, - \, \frac{a_k(k+1)}{2} \, - \, b_k  \\ \vspace{0.2cm}
& = &  \frac{q}{2}  \, - \, \frac{(k-2)(k+1)}{2(k^2+k+2)} \, - \, \frac{k+2}{k^2+k+2}  \\ \vspace{0.2cm}
& = &  \frac{q}{2}  \, - \, \frac{k^2+k+2}{2(k^2+k+2)}  \\ \vspace{0.2cm}
& = & \frac{q-1}{2}.
\end{array}
\]

This is a contradiction which implies that $s \ge 2$.
Recall that $C$ is an odd component in $G-X$ with $s$ edges to $X$. Let $H_{k+1}$ denote $K_{k+1}-e$ and let the two vertices of $H_{k+1}$ with degree $k-1$ be called link vertices.
Let $x_1u_1, x_2u_2, \ldots , x_s u_s$ be distinct edges from $x_i \in X$ to $u_i \in V(C)$ for $i \in [s]$.
Let $G'$ be obtained from $G$ by adding $s-1$ vertex disjoint copies, $G_1,G_2,\ldots,G_{s-1}$, of $H_{k+1}$, and deleting the edges $x_1u_1,x_2u_2, \ldots, x_{s-1}u_{s-1}$ and adding an edge from $x_i$
to a link vertex of $G_i$ for all $i \in [s-1]$. Since no component of $G$ is $k$-regular, it is not difficult to see that no component of $G'$ is $k$-regular.
Note that the following holds.

\begin{itemize}
\item $|V(G')| = |V(G)|+(s-1)(k+1)$
\item $|E(G')| =  |E(G)| + (s-1) \left( \frac{k(k+1)}{2} -1 \right)$
\item $\oc(G'-X) =  \oc(G-X)+s-1$
\end{itemize}


This implies the following.

\[
\begin{array}{rcl}
\multicolumn{3}{l}{\vspace{0.2cm} \frac{1}{2} \left( |V(G')| + |X| - \oc(G' - X)  \right)  } \\ \vspace{0.2cm}
& = &  \frac{1}{2}((|V(G)|+(s-1)(k+1)) + |X| - (\oc(G-X)+s-1)) \\ \vspace{0.2cm}
& = & \frac{1}{2} \left(|V(G)| +|X| - \oc(G-X) \right) + \frac{(s-1)(k+1)}{2} - \frac{s-1}{2} \\ \vspace{0.2cm}
& < & b_k \cdot |E(G)| \, - \, a_k \cdot |V(G)| + \frac{k(s-1)}{2} \\ \vspace{0.2cm}
& = & b_k \cdot \left( |E(G')|-(s-1) \left(\frac{k^2+k-2}{2} \right) \right) \, - \, a_k \cdot \left( |V(G')| - (s-1)(k+1) \right) + \frac{k(s-1)}{2} \\ \vspace{0.2cm}
& = & b_k \cdot |E(G')| \, - \, a_k \cdot |V(G')| \, + \, (s-1) \left( \frac{k}{2} - b_k \left( \frac{k^2+k-2}{2} \right) + a_k (k+1) \right).
\end{array}
\]

As $\xi_i(G',X) = \xi_i(G,X)$ for $i \in [3]$ and $\xi_4(G',X) < \xi_4(G,X)$  we would obtain a contradiction to the lexicographical minimality of
$(\xi_1(G,X), \xi_2(G,X), \ldots, \xi_5(G,X))$ if the following holds.

\[
\begin{array}{crcl}
& 0 & \ge &  \frac{k}{2} - b_k \left( \frac{k^2+k-2}{2} \right) + a_k (k+1) \\
\Updownarrow & & & \\
& \frac{k}{2} & \le &  \frac{k+2}{k^2+k+2} \times \frac{k^2+k-2}{2} - \left( \frac{k-2}{k^2+k+2} \right) (k+1) \\
\Updownarrow & & & \\
& k(k^2+k+2) & \le &  (k+2)(k^2+k-2) - 2(k-2)(k+1) \\
\Updownarrow & & & \\
& k^3+k^2+2k & \le &  k^3 + 3k^2 - 4 - (2k^2-2k-4)  \\
\Updownarrow & & & \\
& k^3+k^2+2k & \le &  k^3 + k^2 + 2k   \\
\end{array}
\]

As the last statement is clearly true, Claim~IV is proved.~\smallqed

\begin{unnumbered}{Claim~V}
$\xi_5(G,X)=0$.
\end{unnumbered}
\proof Suppose, to the contrary, that $\xi_5(G,X) \ge 1$. Let $C$ be a component in $G-X$ which in not isomorphic to $K_{k+1}-e$.
By Claim~II and Claim~IV we note that $C$ is odd and has exactly one edge to $X$. Let $q = |V(C)|$ and let $G'$ be obtained from $G$ by deleting $C$ and adding a copy of $K_{k+1}-e$ and adding an edge from a
degree $k-1$ vertex in $K_{k+1}-e$ to the vertex of $X$ that was adjacent to a vertex of $C$ in $G$.  The following now holds.

\begin{itemize}
\item $|V(G')|=|V(G)| + k + 1 -q $
\item $|E(G')| = |E(G)| + \frac{k(k+1)}{2} - 1 - |E(C)|$
\item $\oc(G'-X) = \oc(G-X)$
\end{itemize}

Since no component of $G$ is $k$-regular, it is not difficult to see that no component of $G'$ is $k$-regular. The following now holds.

\[
\begin{array}{rcl}
\multicolumn{3}{l}{\vspace{0.2cm} \frac{1}{2} \left( |V(G')| + |X| - \oc(G' - X)  \right)  ) } \\ \vspace{0.2cm}
& = & \frac{1}{2}((|V(G)|+ k + 1 - q ) + |X| - (\oc(G-X))
\\ \vspace{0.2cm}
& = & \frac{1}{2} \left(|V(G)| +|X| - \oc(G-X) \right) + \frac{k+1-q}{2} \\ \vspace{0.2cm}
& < & b_k \cdot |E(G)| \, - \, a_k \cdot |V(G)| + \frac{k+1-q}{2}  \\ \vspace{0.2cm}
& = & b_k \cdot \left( |E(G')| - \frac{k(k+1)}{2} + 1 + |E(C)| \right)  \, - \, a_k \cdot \left( |V(G')| - (k +1 -q)  \right) + \frac{k+1-q}{2}  \\ \vspace{0.2cm}
& = & b_k \cdot |E(G')| \, - \, a_k \cdot |V(G')| + \frac{k+1-q}{2} + b_k \left( |E(C)| - \frac{k^2+k-2}{2} \right) + a_k (k+1-q).
\end{array}
\]

We will now consider the cases when $q \le k$ and $q \ge k+1$ seperately. First consider the case when $q \le k$. In this case, $|E(C)| \le q(q-1)/2$ and the above implies the following.

\[
\begin{array}{rcl}
\multicolumn{3}{l}{\vspace{0.2cm} \frac{1}{2} \left( |V(G')| + |X| - \oc(G' - X)  \right)  } \\ \vspace{0.2cm}
& < & b_k \cdot |E(G')| \, - \, a_k \cdot |V(G')| + \frac{k+1-q}{2} + b_k \left( |E(C)| - \frac{k^2+k-2}{2} \right) + a_k (k+1-q)  \\ \vspace{0.2cm}
& \le & b_k \cdot |E(G')| \, - \, a_k \cdot |V(G')| + \frac{k+1-q}{2} + b_k \left( \frac{q(q-1)}{2} - \frac{k^2+k-2}{2} \right) + a_k (k+1-q).
\end{array}
\]

As $\xi_i(G',X) = \xi_i(G,X)$ for $i \in [4]$ and $\xi_5(G',X) < \xi_5(G,X)$  we obtain a contradiction to the lexicographical minimality of
$(\xi_1(G,X), \xi_2(G,X), \ldots, \xi_5(G,X))$ if the following holds

\[
\begin{array}{crcl}
&  \frac{k+1-q}{2} + b_k \left( \frac{q(q-1)}{2} - \frac{k^2+k-2}{2} \right) + a_k (k+1-q)  & \le & 0  \\
\Updownarrow & &  \\
& -q + b_k q (q-1) - 2 a_k q  & \le & b_k(k^2+k-2) - (k+1) - 2a_k(k+1)   \\
\end{array}
\]

Recall that $1 \le q \le k$. Differentiating the left-hand-side twice with respect to $q$ we note that the largest value of $-q + b_k q (q-1) - 2 a_k q $ in the interval $[1,q]$ is obtained at the end points of the interval, when $q=1$ or $q=k$. When $q=1$,
\[
-q + b_k q (q-1) - 2 a_k q = -1-2a_k
\]
and when $q=k$, we get
\[
\begin{array}{lcl}
-q + b_k q (q-1) - 2 a_k q & = & -k + b_k k (k-1) - 2 a_k k \1 \\
& = & -1 + (k-1)(-1 + b_k k - a_k) - a_k(k+1) \1 \\
& = & -1 - a_k(k+1) \1 \\
& < &  -1 - 2a_k.
\end{array}
\]
Therefore the following holds.

\[
\begin{array}{crcl}
& \multicolumn{3}{l}{\vspace{0.2cm}  \frac{k+1-q}{2} + b_k \left( \frac{q(q-1)}{2} - \frac{k^2+k-2}{2} \right) + a_k (k+1-q)   \le  0 } \\
\Updownarrow & &  \\
& -q + b_k q (q-1) - 2 a_k q  & \le & b_k(k^2+k-2) - (k+1) - 2a_k(k+1)   \\
\Uparrow & & & \\
& -1 - 2 a_k  & \le & b_k(k^2+k-2) - (k+1) - 2a_k(k+1)  \\
\Updownarrow & & & \\
& 0 & \le & b_k(k^2+k) -2 b_k - (k+1) - 2a_k(k+1) + 1 + 2a_k \\
\Updownarrow & & & \\
& 0 & \le & (k+1)(k b_k -a_k -1 )  -2 b_k - a_k(k+1) + 1 + 2a_k \\
\Updownarrow & & & \\
& 0 & \le &  -2 b_k - (k-1) a_k + 1  \\
\Updownarrow & & & \\
& 0 & \le &  \frac{-2 (k+2)}{k^2+k+2} - \frac{(k-1)(k-2)}{k^2+k+2}  + \frac{k^2+k+2}{k^2+k+2}  \\
\Updownarrow & & & \\
& 0 & \le &  \frac{(-2k - 4 ) - (k^2-3k+2) + (k^2+k+2) }{k^2+k+2}  \\
\Updownarrow & & & \\
& 0 & \le &  \frac{2k-4}{k^2+k+2}  \\
\end{array}
\]
As the last statement is true we have completed the proof for the case when $q \le k$.
Now consider the case when $q \ge k+1$.
In this case $|E(C)| \le \frac{qk}{2} - 1$ (as $C$ is not $k$-regular) and the following holds.
\[
\begin{array}{rcl}
\multicolumn{3}{l}{\vspace{0.2cm} \frac{1}{2} \left( |V(G')| + |X| - \oc(G' - X)  \right)   } \\ \vspace{0.2cm}
& < &  b_k \cdot |E(G')| \, - \, a_k \cdot |V(G')| + \frac{k+1-q}{2} + b_k \left( |E(C)| - \frac{k^2+k-2}{2} \right) + a_k (k+1-q)  \\ \vspace{0.2cm}
& \le & b_k \cdot |E(G')| \, - \, a_k \cdot |V(G')| + \frac{k+1-q}{2} + b_k \left( \frac{qk-2}{2} - \frac{k^2+k-2}{2} \right) + a_k (k+1-q).
\end{array}
\]

As $\xi_i(G',X) = \xi_i(G,X)$ for $i \in [4]$ and $\xi_5(G',X) < \xi_5(G,X)$  we obtain a contradiction to the lexicographical minimality of
$(\xi_1(G,X), \xi_2(G,X), \ldots, \xi_5(G,X))$ if the following holds

\[
\begin{array}{crcl}
&  \multicolumn{3}{l}{\vspace{0.2cm} \frac{k+1-q}{2} + b_k \left( \frac{qk-2}{2} - \frac{k^2+k-2}{2} \right) + a_k (k+1-q)   \le  0 }  \\
\Updownarrow & &  \\
& q (-1 + kb_k - 2a_k) & \le & b_k(k^2+k-2 + 2) - (k+1) - 2a_k(k+1)   \\
\Updownarrow & &  \\
& q (-1 + 1 -a_k ) & \le & b_k k (k+1) - (k+1) - 2a_k(k+1)   \\
\Uparrow & &  \\
& (k+1) (-a_k) & \le & (k+1) \left( b_k k  - 1 - 2a_k  \right)  \\
\Updownarrow & &  \\
& -a_k & \le &  (b_k k - a_k)  - 1 - a_k  \\
\Updownarrow & &  \\
& -a_k & \le & -a_k  \\
\end{array}
\]
The last statement is clearly true which completes the proof of Claim~V.~\smallqed


\bigskip
We now return to the proof of Theorem~\ref{match_thm2}. As mentioned before the statement of Claim~I, we may assume that $G$ is connected. By Claim~I, $\xi_1(G,X)=0$, and so there is no edge in $G[X]$. By Claim~IV, $\xi_4(G,X)=0$, and so every component in $C-X$ has exactly one edge to $X$. Since $C$ is connected, we therefore note that $|X|=1$. Let $X=\{x\}$. By Claim~III, $\xi_3(G,X)=0$, and so $d_G(x)=k$. By Claim~V, every component in $C-X$ is isomorphic to $K_{k+1}-e$, which implies the following.

\begin{itemize}
\item $|V(G)| = k (k+1) + 1 = k^2+k+1$.
\item $|E(G)| = k \left( \frac{k(k+1)}{2} -1 \right) + d(x) = k \left( \frac{k^2+k-2}{2} \right) + k$.
\item $\alpha'(G) = k \left( \frac{k}{2} \right) + 1 = \frac{k^2+2}{2}$.
\end{itemize}

Therefore, the following holds.

\[
\begin{array}{rcl} \vspace{0.2cm}
b_k |E(G)|- a_k |V(G)|
& = & \frac{k+2}{k^2+k+2} \times k \times \left( \frac{k^2+k-2}{2} + 1 \right) - \frac{k-2}{k^2+k+2} \times (k^2+k+1) \\ \vspace{0.2cm}
& = & \frac{1}{2(k^2+k+2)} \left( k(k+2)(k^2+k) - 2(k-2)(k^2+k+2) \right) \\ \vspace{0.2cm}
& = & \frac{ (k^4+3k^3+2k^2)-(2k^3 -2k^2 -2k -4)}{2(k^2+k+2)} \\ \vspace{0.2cm}
& = & \frac{ k^4+k^3+4k^2+2k+4}{2(k^2+k+2)} \\ \vspace{0.2cm}
& = & \frac{ (k^2+2)(k^2+k+2)}{2(k^2+k+2)} \\ \vspace{0.2cm}
& = & \frac{ k^2+2}{2} \\ \vspace{0.2cm}
& = & \alpha'(G).
\end{array}
\]

This contradicts the fact that $G$ was a counter example to the theorem.~\qed

\section{Proof of Corollary~\ref{corKeven}}
\label{S:corKeven}

Recall the statement of Corollary~\ref{corKeven}.

\noindent \textbf{Corollary~\ref{corKeven}} \emph{If $k \ge 4$ is an even integer and $G$ is a connected graph of order $n$, size $m$ and maximum degree $\Delta(G) \le k$, then
\[
\alpha'(G) \ge \frac{n}{k(k+1)} \, + \, \frac{m}{k+1} - \frac{1}{k(k+1)},
\]
unless the following holds.
\\[-24pt]
\begin{enumerate}
\item $G$ is $k$-regular and $n = k+1$, in which case $\alpha'(G) = \frac{n-1}{2} = \frac{n}{k(k+1)} \, + \, \frac{m}{k+1} - \frac{1}{k}$.
\item $G$ is $k$-regular and $n = k+3$, in which case $\alpha'(G) = \frac{n-1}{2} = \frac{n}{k(k+1)} \, + \, \frac{m}{k+1} - \frac{3}{k(k+1)}$.
\end{enumerate}
}

\noindent \textbf{Proof of Corollary~\ref{corKeven}.}  If $G$ is not $k$-regular, then the result follows from Theorem~\ref{match_thm}, so assume that $G$ is $k$-regular. By the $k$-regularity of $G$ we have $nk = 2m$, which together with the observation that
\[
\frac{k^2+2}{k^2+k} \le \frac{k^2+4}{k^2+k+2}.
\]
when $k \ge 4$, implies the following.

\[
\begin{array}{rcl} \vspace{0.2cm}
\displaystyle{ \frac{n}{k(k+1)} \, + \, \frac{m}{k+1} - \frac{1}{k(k+1)} }
    & = &  \displaystyle{   \frac{n}{k(k+1)} \, + \, \frac{nk}{2(k+1)} - \frac{1}{k(k+1)} } \\ \vspace{0.2cm}
    & = &  \displaystyle{ \frac{n}{2} \left(  \frac{2}{k(k+1)} \, + \, \frac{k}{k+1} \right) - \frac{1}{k(k+1)} } \\ \vspace{0.2cm}
    & < &  \displaystyle{ \frac{n}{2} \left(  \frac{2}{k(k+1)} \, + \, \frac{k^2}{k(k+1)} \right) } \\ \vspace{0.2cm}
    & = &  \displaystyle{ \frac{n}{2} \times  \frac{k^2+2}{k(k+1)}.} \\ \vspace{0.2cm}
    & \le &  \displaystyle{ \frac{n}{2} \times  \frac{k^2+4}{k^2+k+2}. }
\end{array}
\]

By Theorem~\ref{KregEVEN} we therefore have the following.

\[
\alpha'(G) \ge \min \left\{ \frac{k^2+4}{k^2+k+2} \times \frac{n}{2} , \frac{n-1}{2} \right\} \geq
\min \left\{ \frac{n}{k(k+1)} \, + \, \frac{m}{k+1} - \frac{1}{k(k+1)}  , \frac{n-1}{2} \right\}.
\]

As $\frac{n}{k(k+1)} + \frac{m}{k+1} - \frac{1}{k(k+1)} < \frac{n}{2}$ by the above, this proves the theorem when $n$ is even.
So let $n$ be odd. We will now determine when the following holds.

\[
\frac{n}{k(k+1)} \, + \, \frac{nk}{2(k+1)} - \frac{1}{k(k+1)} > \frac{n-1}{2}.
\]

Define $r$ such that $n = k + r$ and note that the above is equivalent to the following.

\[
\begin{array}{rcl} \vspace{0.2cm}
0 & < &  \displaystyle{ \frac{n}{k(k+1)} \, + \, \frac{nk}{2(k+1)} - \frac{1}{k(k+1)} - \frac{n-1}{2} } \\ \vspace{0.2cm}
& = & \displaystyle{ n \left(\frac{ 2 + k^2 - k(k+1)}{2k(k+1)}\right) - \frac{1}{k(k+1)} + \frac{1}{2} } \\ \vspace{0.2cm}
& = & \displaystyle{ (k+r) \left( \frac{2-k}{2k(k+1)} \right) + \frac{-2+k(k+1)}{2k(k+1)} } \\ \vspace{0.2cm}
& = & \displaystyle{ \frac{-k^2+2k+2r-kr}{2k(k+1)} + \frac{k^2+k-2}{2k(k+1)} } \\ \vspace{0.2cm}
& = & \displaystyle{ \frac{3k+2r-kr-2}{2k(k+1)}. }
\end{array}
\]

This is equivalent to $3k-2-r(k-2) > 0$, which is equivalent to the following (as $k \ge 4$).

\[
  r < \frac{3k-2}{k-2}.
\]

When $k \ge 4$ this is equivalent to $r < 5$ (as $r$ is odd) and therefore the following holds (as we already handled the case when $n$ was even).

\[
\begin{array}{rcl} \vspace{0.4cm}
\alpha'(G) & \ge  &  \displaystyle{ \min \left\{ \frac{n}{k(k+1)} \, + \, \frac{m}{k+1} - \frac{1}{k(k+1)}  , \frac{n-1}{2} \right\} } \\
& = &
\left\{
\begin{array}{ccl} \vspace{0.4cm}
\displaystyle{ \frac{n}{k(k+1)} \, + \, \frac{m}{k+1} - \frac{1}{k(k+1)} } & \hspace{0.3cm} & \mbox{if } n \not\in \{k+1,k+3\} \\
\displaystyle{ \frac{n-1}{2} } & &  \mbox{if } n \in \{k+1,k+3\}.  \\
\end{array}
\right. \\
\end{array}
\]

We will therefore determine $\theta$ such that the following holds.

\[
\begin{array}{crcl} \vspace{0.2cm}
& \displaystyle{  \frac{n-1}{2} } & = & \displaystyle{ \frac{n}{k(k+1)} \, + \, \frac{m}{k+1} - \frac{\theta}{k(k+1)} } \\ \vspace{0.2cm}
\Updownarrow & & & \\ \vspace{0.2cm}
& 2 \theta & = & 2n  + 2km - (n-1)(k(k+1)) \\ \vspace{0.2cm}
& & = & 2n  + nk^2 - n(k^2+k) + (k^2+k) \\ \vspace{0.2cm}
& & = & n (2-k) + k^2+k \\ \vspace{0.2cm}
& & = & (k+r) (2-k) + k^2+k \\ \vspace{0.2cm}
& & = & 3k +2r -rk.
\end{array}
\]
So, if $r =3$, then $\theta=3$ and if $r=1$, then $\theta = k+1$. This completes the proof of Corollary~\ref{corKeven}.~\qed

\medskip
The lower bound in Corollary~\ref{corKeven} is tight when $G$ is $k$-regular and $n \in \{k+1,k+3\}$. We will furthermore show that they are tight for infinite classes of trees and infinite classes
of graphs with any given average degree between $2$ and $k-\frac{k-2}{k^2}$, implying that the corollaries are tight for almost all possible densities. Note that for small values of $k$ the value $k-\frac{k-2}{k^2}$ is the following.

\begin{center}
\[
\begin{array}{|c||c|c|c|c|c|c|} \hline
k &                  4     &     6         &     8         &  10   &    12         & 14             \\ \hline \hline
k-\frac{k-2}{k^2} & 3.875  & 5.888 \ldots  & 7.906 \ldots  & 9.92  & 11.93 \ldots  & 13.938 \ldots  \\ \hline
\end{array}
\]
\end{center}
\begin{center}
\textbf{Table~3.} The value $k-\frac{k-2}{k^2}$ for small values of $k \ge 4$ with $k$ even.
\end{center}

We will illustrate how to obtain the above.
Let $k \ge 4$ be even and let $r \ge 1$ be arbitrary and let $\ell = r(k-1)+1$. Let $X_1,X_2,\ldots,X_\ell$ be a number of vertex disjoint graphs such that each $X_i$ where $i \in [\ell]$ is either a
single vertex or it is a $K_{k+1}$ where an arbitrary edge has been deleted.
Let $Y=\{y_1,y_2,\ldots,y_r\}$ and build the graph $G_{k,r}$ as follows. Let $G_{k,r}$ be obtained from the disjoint union of the graphs  $X_1,X_2,\ldots,X_\ell$ by adding to it the vertices in $Y$ and furthermore, for every $i \in [r]$, adding an edge from $y_i$ to a vertex in each graph $X_{(i-1)(k-1)+1}$, $X_{(i-1)(k-1)+2}$, $X_{(i-1)(k-1)+3}, \ldots, X_{(i-1)(k-1)+k}$ in such a way that no vertex degree becomes more than~$k$. Let $\cG_{k,r}$ be the family of all such graph $G_{k,r}$.

When $k=4$ and $r=2$, an example of a graph $G$ in the family $\cG_{k,r}$ is illustrated in Figure~\ref{graph1}, where $G$ has order~$n = 21$, size~$m = 35$ and matching number $\alpha'(G)=8$.

\begin{figure}[htb]
\begin{center}
\tikzstyle{vertexX}=[circle,draw, fill=black!90, minimum size=8pt, scale=0.5, inner sep=0.2pt]
\begin{tikzpicture}[scale=0.21]
 \node (a1) at (1.0,3.0) [vertexX] {};
\node (a2) at (2.0,6.0) [vertexX] {};
\node (a3) at (5.0,6.0) [vertexX] {};
\node (a4) at (6.0,3.0) [vertexX] {};
\node (a5) at (3.5,1.0) [vertexX] {};
\node (b1) at (9.0,4.0) [vertexX] {};
\node (c1) at (12.0,4.0) [vertexX] {};
\node (d1) at (15.0,3.0) [vertexX] {};
\node (d2) at (16.0,6.0) [vertexX] {};
\node (d3) at (19.0,6.0) [vertexX] {};
\node (d4) at (20.0,3.0) [vertexX] {};
\node (d5) at (17.5,1.0) [vertexX] {};
\node (e1) at (23.0,4.0) [vertexX] {};
\node (f1) at (26.0,3.0) [vertexX] {};
\node (f2) at (27.0,6.0) [vertexX] {};
\node (f3) at (30.0,6.0) [vertexX] {};
\node (f4) at (31.0,3.0) [vertexX] {};
\node (f5) at (28.5,1.0) [vertexX] {};
\node (g1) at (34.0,4.0) [vertexX] {};
\node (t1) at (10.5,12.0) [vertexX] {};
\node (t2) at (28.0,12.0) [vertexX] {};
\draw (a2) -- (a1) -- (a5) -- (a4) -- (a3) -- (a1) -- (a4) -- (a2) -- (a5) -- (a3);
\draw (d2) -- (d1) -- (d5) -- (d4) -- (d3) -- (d1) -- (d4) -- (d2) -- (d5) -- (d3);
\draw (f2) -- (f1) -- (f5) -- (f4) -- (f3) -- (f1) -- (f4) -- (f2) -- (f5) -- (f3);
\draw (a3) -- (t1) -- (b1);
 \draw (c1) -- (t1) -- (d2);
\draw (g1) -- (t2) -- (f3);
 \draw (e1) -- (t2) -- (d3);
 \end{tikzpicture}
\end{center}
\vskip -0.5 cm \caption{A graph $G$ in the family $\cG_{4,2}$}  \label{graph1}
\end{figure}

\begin{prop}
For $k \ge 4$ an even integer and $r \ge 1$ arbitrary, if  $G \in \cG_{k,r}$ has order~$n$ and size~$m$, then \[
\alpha'(G) = \frac{n}{k(k+1)} \, + \, \frac{m}{k+1} - \frac{1}{k}.
\]
\label{p:Gkr}
\end{prop}
\proof Assume that in $X_1,X_2,\ldots,X_\ell$ we have $\ell_1$ single vertices and $\ell_2$ copies of $K_{k+1}$'s minus an edge. Note that $\ell = \ell_1 + \ell_2$ and
$n = r + \ell_1 + \ell_2(k+1)$. Furthermore we have $m = rk + \ell_2(k(k+1)/2 -1)$ and $\alpha'(G) = r + \ell_2(k/2)$. Therefore the following holds.

\[
\begin{array}{rcl}
\multicolumn{3}{l}{\vspace{0.2cm} \displaystyle{ \frac{n}{k(k+1)} \, + \, \frac{m}{k+1} - \frac{1}{k(k+1)} }} \\ \vspace{0.2cm}
& = & \displaystyle{ \frac{r + \ell_1 + \ell_2(k+1) }{k(k+1)} \, + \, \frac{2rk + \ell_2(k^2+k-2)}{2(k+1)} - \frac{1}{k(k+1)} } \\ \vspace{0.2cm}
& = & \displaystyle{ \frac{r + (\ell-\ell_2) + \ell_2k + \ell_2 }{k(k+1)} \, + \, \frac{2rk^2 + \ell_2k(k^2+k-2)}{2k(k+1)} - \frac{2}{2k(k+1)} } \\ \vspace{0.2cm}
& = & \displaystyle{ \frac{2r + 2\ell + 2\ell_2k   + 2rk^2 + \ell_2k(k^2+k-2) - 2}{2k(k+1)} } \\ \vspace{0.2cm}
& = & \displaystyle{ \frac{2r + 2(r(k-1)+1) +  2rk^2 + \ell_2k(k^2+k) - 2}{2k(k+1)} } \\ \vspace{0.2cm}
& = & \displaystyle{ \frac{2r + 2rk-2r+2 +  2rk^2 + \ell_2k^2(k+1) - 2}{2k(k+1)} } \\ \vspace{0.2cm}
& = & \displaystyle{ \frac{ 2r(k^2+k) + \ell_2k^2(k+1) }{2k(k+1)} } \\ \vspace{0.2cm}
& = & \displaystyle{ r +  \frac{1}{2} k \ell_2} \\
& = & \alpha'(G). \hspace*{0.5cm} \Box\\
\end{array}
\]

By Proposition~\ref{p:Gkr}, the lower bound on the matching number in Corollary~\ref{corKeven} is tight for every graph in the family~$\cG_{k,r}$. We remark that if all $X_1,X_2,\ldots,X_\ell$ used to construct a graph in the family~$\cG_{k,r}$ are single vertices, then clearly we have a tree and as $r \ge 1$ was arbitrary we obtain an infinite class of trees where the corollaries are tight. If all $X_1,X_2,\ldots,X_\ell$ are copies of $K_{k+1}$ minus an edge, then we denote the resulting subfamily of graphs of $\cG_{k,r}$ by $\cG_{k,r}'$. Hence, each graph in the family $\cG_{k,r}'$ achieves the lower bound in Corollary~\ref{corKeven}.
We remark, further, that if we build a graph in the family~$\cG_{k,r}$ using only single vertices for the copies of $X_i$ for each $i \in [\ell]$, then the resulting graph $G$ is a tree. Hence, as an immediate consequence of Proposition~\ref{p:Gkr} and the above observations, we have the following result.

\begin{prop}
The lower bound in Corollary~\ref{corKeven} is achieved for both trees and for the class of graphs in the family~$\cG_{k,r}'$.
\label{p:Gkrtree}
\end{prop}

We note that the average degree of a graph in the family $\cG_{k,r}'$ is the following, as $n=\ell(k+1)+r$ and there are $\ell-(r-1)$ vertices of degree $k-1$.
\[
\begin{array}{rcl}
\displaystyle{ \frac{1}{n} \sum_{v \in V(G)} d(v) } & = & \displaystyle{\frac{1}{n} \left( nk - (\ell-(r-1)) \right) } \\ \vspace{0.2cm}
& = & \displaystyle{k - \frac{1}{n}(\ell-r+1) } \\ \vspace{0.2cm}
& = & \displaystyle{k - \frac{r(k-1)+1-r+1}{(r(k-1)+1)(k+1)+r} }\\ \vspace{0.2cm}
& = & \displaystyle{k - \frac{r(k-2)+2}{rk^2 + k+1}.} \\
\end{array}
\]

So when $r$ is large the average degree can get arbitrarily close to $k - \frac{k-2}{k^2}$. Clearly, as the average degree of a tree is less than~$2$, the average degree can be made arbitrarily close to any $\beta$ where $2 \le \beta < k-\frac{k-2}{k^2}$ by picking the correct proportion of single vertices in $X_1,X_2,\ldots,X_\ell$.

\section{Proof of Corollary~\ref{corKodd}}
\label{S:corKodd}

Recall the statement of Corollary~\ref{corKodd}.

\medskip
\noindent \textbf{Corollary~\ref{corKodd}} \emph{If $k \ge 3$ is an odd integer and $G$ is a connected graph of order $n$, size $m$, and with maximum degree $\Delta(G) \le k$, then
\[
\alpha'(G) \ge \left( \frac{k-1}{k(k^2 - 3)} \right) n \, + \, \left( \frac{k^2 - k - 2}{k(k^2 - 3)} \right) m \, - \, \frac{k-1}{k(k^2 - 3)}.
\]
}

\noindent \textbf{Proof of Corollary~\ref{corKodd}.}  If $G$ is not $k$-regular, then the result follows from Theorem~\ref{match_thm}, so assume that $G$ is $k$-regular. By the $k$-regularity of $G$ we have $nk = 2m$, which implies the following by Theorem~\ref{KregODD}.

\[
\begin{array}{rcl}
\multicolumn{3}{l}{\vspace{0.2cm} \displaystyle{ \left( \frac{k-1}{k(k^2 - 3)} \right) n \, + \, \left( \frac{k^2 - k - 2}{k(k^2 - 3)} \right) m \, - \, \frac{k-1}{k(k^2 - 3)} }} \3 \\
& = &  \displaystyle{   \left( \frac{k-1}{k(k^2 - 3)} \right) n \, + \, \left( \frac{k(k^2 - k - 2)}{2k(k^2 - 3)} \right) n \, - \, \frac{k-1}{k(k^2 - 3)} } \3 \\
& = &  \displaystyle{ \frac{(k^3-k^2-2) \, n - 2k + 2}{2k(k^2-3)} } \1 \\
& \stackrel{(Thm~\ref{KregODD})}{\le}  &  \alpha'(G). \hspace*{0.5cm} \Box
\end{array}
\]

We show that the lower bound on the matching number in Corollary~\ref{corKodd} is tight for infinite classes of trees and other infinite classes (including the class of $k$-regular graphs) of connected graphs with maximum degree at most~$k$.

For $k \ge 3$ odd, let $H_{k+2}$ be the graph of (odd)  order~$k+2$ obtained from $K_{k+2}$ by removing the edges of an almost perfect matching; that is, the complement $\overline{H_{k+2}}$ of $H_{k+2}$ is isomorphic to $P_3 \cup  ( \frac{k-1}{2} )P_2$. We note that every vertex in $H_{k+2}$ has degree~$k$, except for exactly one vertex, which has degree~$k-1$. We call the vertex of degree~$k-1$ in $H_{k+2}$ the \emph{link vertex} of $H_{k+2}$. We note that $H_{k+2}$ has size~$m(H_{k+2}) = \frac{1}{2}(k^2 + 2k - 1)$.

For $k \ge 3$ odd and $r \ge 1$ arbitrary, let $T_{k,r}$ be a tree with maximum degree at most~$k$ and with partite sets $V_1$ and $V_2$, where $|V_2| = r$. Let $H_{k,r}$ be obtained from $T_{k,r}$ as follows: For every vertex $x$ in $V_2$ with $d_{T_{k,r}}(x) < k$, add $k-d_{T_{k,r}}(x)$ copies of the subgraph $H_{k+2}$ to $T_{k,r}$ and in each added copy of $H_{k+2}$, join the link vertex of $H_{k+2}$ to $x$. We note that every vertex in the resulting graph $H_{k,r}$ has degree~$k$, except possibly for vertices in the set $V_1$ whose degrees belong to the set $\{1,2,\ldots,k\}$. Let $\cH_{k,r}$ be the family of all such graph $H_{k,r}$.

When $k=3$ and $r=4$, an example of a graph $G$ in the family $\cH_{k,r}$ is illustrated in Figure~\ref{H34}, where $G$ has order~$n = 29$, size~$m = 40$ and matching number $\alpha'(G)=12$.

\begin{figure}[htb]
\begin{center}
\tikzstyle{vertexX}=[circle,draw, fill=black!90, minimum size=8pt, scale=0.5, inner sep=0.2pt]
\begin{tikzpicture}[scale=0.21]
\node (a0) at (3,7.5) [vertexX] {};
\node (a1) at (1,3.0) [vertexX] {};
\node (a2) at (5,3.0) [vertexX] {};
\node (a3) at (3,5.0) [vertexX] {};
\node (a4) at (3,1.0) [vertexX] {};
\draw (a3) -- (a4) -- (a1) -- (a0) -- (a2) -- (a4);
\draw (a1) -- (a3) -- (a2);
\node (b0) at (9,7.5) [vertexX] {};
\node (b1) at (7,3.0) [vertexX] {};
\node (b2) at (11,3.0) [vertexX] {};
\node (b3) at (9,5.0) [vertexX] {};
\node (b4) at (9,1.0) [vertexX] {};
\draw (b3) -- (b4) -- (b1) -- (b0) -- (b2) -- (b4);
\draw (b1) -- (b3) -- (b2);
\node (c0) at (15,7.5) [vertexX] {};
\node (c1) at (13,3.0) [vertexX] {};
\node (c2) at (17,3.0) [vertexX] {};
\node (c3) at (15,5.0) [vertexX] {};
\node (c4) at (15,1.0) [vertexX] {};
\draw (c3) -- (c4) -- (c1) -- (c0) -- (c2) -- (c4);
\draw (c1) -- (c3) -- (c2);
\node (d0) at (27,7.5) [vertexX] {};
\node (d1) at (25,3.0) [vertexX] {};
\node (d2) at (29,3.0) [vertexX] {};
\node (d3) at (27,5.0) [vertexX] {};
\node (d4) at (27,1.0) [vertexX] {};
\draw (d3) -- (d4) -- (d1) -- (d0) -- (d2) -- (d4);
\draw (d1) -- (d3) -- (d2);
\node (x0) at (2.0,17) [vertexX] {};
\node (x1) at (8.5,17) [vertexX] {};
\node (x2) at (15.0,17) [vertexX] {};
\node (x3) at (21.5,17) [vertexX] {};
\node (x4) at (28.0,17) [vertexX] {};
\node (y0) at (3,12) [vertexX] {};
\node (y1) at (11,12) [vertexX] {};
\node (y2) at (19,12) [vertexX] {};
\node (y3) at (27,12) [vertexX] {};
\draw (y0) -- (a0);
\draw (y1) -- (b0);
\draw (y1) -- (c0);
\draw (y3) -- (d0);
\draw (x0) -- (y0) -- (x1) -- (y1);
\draw (x1) -- (y2) -- (x2);
\draw (y2) -- (x3) -- (y3) -- (x4);

 \draw (31.5,17) node {$V_1$};
 \draw (31,12) node {$V_2$};
 \draw [rounded corners] (0.5,16) rectangle (29.5,18);
 \draw [rounded corners] (1,11) rectangle (29,13);

\end{tikzpicture}
\end{center}
\vskip -0.5 cm \caption{A graph $G$ in the family $\cH_{3,4}$}  \label{H34}
\end{figure}

\begin{prop}
For $k \ge 3$ an odd integer and $r \ge 1$ arbitrary, if  $G \in \cH_{k,r}$ has order~$n$ and size~$m$, then
\[
\alpha'(G) = \left( \frac{k-1}{k(k^2 - 3)} \right) n \, + \, \left( \frac{k^2 - k - 2}{k(k^2 - 3)} \right) m \, - \, \frac{k-1}{k(k^2 - 3)}.
\]
\label{p:Hkr}
\end{prop}
\proof Let $G \cong H_{k,r} \in \cH_{k,r}$ have order~$n$ and size~$m$. Suppose that $\ell$ copies of the graph $H_{k+2}$ were added when constructing the graph $G$. Thus,
\[
\begin{array}{rcl}
\ell & = & \displaystyle{ k|V_2| - \sum_{x \in V_2} d_{T_{k,r}}(x) } \1 \\
& = & k|V_2| - m(T_{k,r}) \1 \\
& = &  k|V_2| - (|V_1| + |V_2| - 1) \1 \\
& = &  (k-1)|V_2| - |V_1| + 1.
\end{array}
\]
The graph $G$ has order
\[
\begin{array}{rcl}
n = n(G) & = & |V_1| + |V_2| + \ell(k+2) \1 \\
& = & |V_1| + |V_2| + ((k-1)|V_2| - |V_1| + 1)(k+2)  \1 \\
& = &  (k^2 + k - 1)|V_2| - (k+1)|V_1| + (k+2).
\end{array}
\]
and size
\[
\begin{array}{rcl}
m = m(G) & = & k|V_2| +  \frac{1}{2} \ell (k^2 + 2k - 1) \1 \\
& = & k|V_2| +  \frac{1}{2} ((k-1)|V_2| - |V_1| + 1) (k^2 + 2k - 1)  \1 \\
& = &  \frac{1}{2}(k^3 + k^2 - k + 1)|V_2| - \frac{1}{2}(k^2 + 2k - 1) |V_1| + \frac{1}{2}(k^2 + 2k - 1).
\end{array}
\]

Furthermore by deleting the vertices $V_2$ from $G$ we obtain
$\ell + |V_1| = (k-1)|V_2| + 1$ odd components. Therefore, by
Theorem~\ref{Berge},
\[
\begin{array}{rcl}
\2
2\alpha'(G) & \le & |V(G)|+|V_2|-\oc(G - V_2)  \\
\2
  & = &  ((k^2 + k - 1)|V_2| - (k+1)|V_1| + (k+2)) + |V_2|  - ((k-1)|V_2| + 1)   \\
\2
  & = & (k^2 + 1)|V_2| - (k+1)|V_1| + (k + 1).
\end{array}
\]
However, the lower bound of Corollary~\ref{corKodd} shows that
\[
\begin{array}{cl}
& 2\alpha'(G) \3 \\
\ge & \displaystyle{ 2\left( \frac{k-1}{k(k^2 - 3)} \right) n \, + \, 2\left( \frac{k^2 - k - 2}{k(k^2 - 3)} \right) m \, - \, \frac{2(k-1)}{k(k^2 - 3)} } \3 \\
= & \displaystyle{ 2\left( \frac{k-1}{k(k^2 - 3)} \right) ((k^2 + k - 1)|V_2| - (k+1)|V_1| + (k+2)) }
\2 \\
 & \hspace*{0.15cm} \displaystyle{ \, + \, \left( \frac{k^2 - k - 2}{k(k^2 - 3)} \right) \left( (k^3 + k^2 - k + 1)|V_2| - (k^2 + 2k - 1) |V_1| + (k^2 + 2k - 1) \right)  } \2 \\
 & \hspace*{0.75cm} \displaystyle{ \, - \, \frac{2(k-1)}{k(k^2 - 3)} } \3 \\
= &  \displaystyle{ \left( \frac{k^5 - 2k^3 - 3k}{k(k^2 - 3)} \right)|V_2| - \left( \frac{k^4 + k^3 - 3k^2 - 3k}{k(k^2 - 3)} \right) |V_1| + \left( \frac{k^4 + k^3 - 3k^2 - 3k}{k(k^2 - 3)} \right) } \3 \\
= &  (k^2 + 1)|V_2| - (k+1)|V_1| + (k + 1).
\end{array}
\]

\noindent
Consequently, we must have equality throughout the above inequality chains. In particular,
\[
\alpha'(G) = \left( \frac{k-1}{k(k^2 - 3)} \right) n \, + \, \left( \frac{k^2 - k - 2}{k(k^2 - 3)} \right) m \, - \, \frac{k-1}{k(k^2 - 3)}.
\]
This completes the proof of the proposition.~\qed

\medskip
We remark that if the tree $T_{k,r}$ used to construct the graph $H_{k,r}$ is chosen so that every vertex in $V_2$ has degree~$k$, then $H_{k,r}$ is a tree. If, however, tree $T_{k,r}$ used to construct the graph $H_{k,r}$ is chosen so that every vertex in $V_1$ has degree~$k$, then $H_{k,r}$ is a $k$-regular graph. Hence, as an immediate consequence of Proposition~\ref{p:Hkr}, we have the following result.

\begin{prop}
The lower bound in Corollary~\ref{corKodd} is achieved for an infinite class of both trees and $k$-regular graphs.
\label{p:Hkrtreeregular}
\end{prop}

\section{Proof of Corollary~\ref{corKeven2}}
\label{S:corKeven2}

Recall the statement of Corollary~\ref{corKeven2}.

\medskip
\noindent \textbf{Corollary~\ref{corKeven2}}
\emph{
If $k \ge 4$ is an even integer and $G$ is a graph of order $n$, size $m$ and maximum degree $\Delta(G) \le k$, then
\[
 \alpha'(G) \ge  \left( \frac{k+2}{k^2+k+2} \right) m  \, - \,  \left( \frac{k-2}{k^2+k+2} \right) n
\]
unless the following holds.
\\[-20pt]
\begin{enumerate}
\item $G$ is $k$-regular and $n = k+1$, in which case  $\alpha'(G) \ge  \left( \frac{k+2}{k^2+k+2} \right) m  \, - \,  \left( \frac{k-2}{k^2+k+2} n \right)  \, - \frac{k+2}{k^2+k+2}$.
\item $G$ is $k$-regular and $n = k+3$, in which case $\alpha'(G) \ge  \left( \frac{k+2}{k^2+k+2} \right) m  \, - \,  \left( \frac{k-2}{k^2+k+2} \right) n  \, - \frac{4}{k^2+k+2}$.
\item $G$ is $4$-regular and $n = 9$, in which case $\alpha'(G) \ge  \left( \frac{k+2}{k^2+k+2} \right) m  \, - \,  \left( \frac{k-2}{k^2+k+2} \right) n  \, - \frac{2}{k^2+k+2}$.
\end{enumerate}
}

\noindent \textbf{Proof of Corollary~\ref{corKeven2}.}
 If $G$ is not $k$-regular, then the result follows from Theorem~\ref{match_thm2}, so assume that $G$ is $k$-regular.
By the $k$-regularity of $G$ we have $nk = 2m$, which implies the following

\[
\begin{array}{rcl} \vspace{0.2cm}
\displaystyle{ \left( \frac{k+2}{k^2+k+2} \right) m \, - \, \left( \frac{k-2}{k^2+k+2} \right) n }
    & = &  \displaystyle{ \frac{k+2}{k^2+k+2} \times \frac{nk}{2} \, - \, \left( \frac{2k-4}{2(k^2+k+2)} \right) n } \\ \vspace{0.2cm}
    & = &  \displaystyle{ \frac{n}{2} \times \frac{k(k+2) - (2k-4)}{k^2+k+2} } \\ \vspace{0.2cm}
    & = &  \displaystyle{ \frac{n}{2} \times \frac{k^2+4}{k^2+k+2} }
\end{array}
\]

Applying Theorem~\ref{KregEVEN} to the $k$-regular graph $G$, the matching number of $G$ is bounded below as follows.
\[
\alpha'(G) \ge \min \left\{ \frac{k^2+4}{k^2+k+2} \times \frac{n}{2} , \frac{n-1}{2} \right\}
\]

Therefore we have shown that the corollary holds in all cases, except when

\[
\frac{k^2+4}{k^2+k+2} \times \frac{n}{2} > \alpha'(G) \ge \frac{n-1}{2}.
\]

So assume that this exceptional case occurs. Since $k \ge 4$, we note that in this case

\[
\frac{n-1}{2} \le \alpha'(G) < \frac{k^2+4}{k^2+k+2} \times \frac{n}{2} < 1 \times \frac{n}{2} = \frac{n}{2},
\]

\noindent
implying that  $\alpha'(G) = (n-1)/2$, and so $n \ge k + 1$ is odd. Thus, $n = k+i$ for some $i \ge 1$ odd. This implies the following.

\[
\begin{array}{rcl} \vspace{0.2cm}
b_k |E(G)|- a_k |V(G)| - \alpha'(G)
& = & \displaystyle{ \left( \frac{k+2}{k^2+k+2} \right) m \, - \, \left( \frac{k-2}{k^2+k+2} \right) n - \alpha'(G) } \\ \vspace{0.2cm}
& = & \displaystyle{ \frac{k^2+4}{k^2+k+2} \times \frac{n}{2} - \frac{n-1}{2} }\\ \vspace{0.2cm}
& = & \displaystyle{ \frac{k^2+4}{k^2+k+2} \times \frac{k+i}{2} - \frac{(k+i-1)(k^2+k+2)}{2(k^2+k+2)} }\\ \vspace{0.2cm}
& = & \displaystyle{\frac{1}{2(k^2+k+2)} \left( (k+i)(k^2+4) - (k+i-1)(k^2+k+2) \right) }\\ \vspace{0.2cm}
& = & \displaystyle{ \frac{1}{2(k^2+k+2)} \left( (k+i)(2-k) + k^2+k+2  \right) }\\ \vspace{0.2cm}
& = & \displaystyle{ \frac{1}{2(k^2+k+2)} \left(  3k+2 - i(k-2) \right),}
\end{array}
\]
or, equivalently,
\begin{equation}
\alpha'(G) = b_k |E(G)|- a_k |V(G)| - \frac{3k+2-i(k-2)}{2(k^2+k+2)}.
\label{Eq1c}
\end{equation}

This proves the case when $n=k+1$, as when $i=1$  the right-hand side of Equation~(\ref{Eq1c}) is equivalent to the lower bound on $\alpha'(G)$  in the statement of Corollary~\ref{corKeven2}(a).
We also get the case when $n=k+3$, as when $i=3$ we again note that the the right-hand side of Equation~(\ref{Eq1c}) is equivalent to the lower bound on $\alpha'(G)$ in the statement of Corollary~\ref{corKeven2}(b).
When $n=k+5$, we have $i=5$ and in this case $3k+2-i(k-2) = 12-2k$. When $k=4$, the right-hand side of Equation~(\ref{Eq1c}) is equivalent to the lower bound on $\alpha'(G)$  in the statement of Corollary~\ref{corKeven2}(c). When $k \ge 6$, we note that $12-2k \le 0$, implying that
\[
\alpha'(G) \ge b_k |E(G)|- a_k |V(G)| = \frac{k^2+4}{k^2+k+2} \times \frac{n}{2},
\]
a contradiction to our exceptional case.
When $n \ge k+7$, then $i \ge 7$ and $3k+2-i(k-2) \le 16-4k \le 0$, implying that
$\alpha'(G) \ge b_k |E(G)|- a_k |V(G)|$, again a contradiction to our exceptional case.~\qed

\medskip
Following the notation in the proof of Proposition~\ref{p:Gkr}, for $k \ge 4$ an even integer and $r \ge 1$ arbitrary, let $G$ be a graph in the family $\cG_{k,r}'$ of order~$n$ and size~$m$. We recall that contains $\ell$ copies of $K_{k+1} - e$ with $\ell = r(k-1)+1$. Note that in this case, $n = r + \ell(k+1)$, $m = rk + \frac{1}{2}\ell(k^2 + k - 2)$ and $\alpha'(G) = r + \frac{1}{2}\ell k$. Therefore the following holds.

\[
\begin{array}{rcl}
\multicolumn{3}{l}{\vspace{0.2cm} \displaystyle{ \left( \frac{k+2}{k^2+k+2} \right) m  \, - \,  \left( \frac{k-2}{k^2+k+2} \right) n  }} \\ \vspace{0.2cm}
& = & \displaystyle{ \left( \frac{k+2}{k^2+k+2} \right) \left( rk + \frac{\ell(k^2 + k - 2)}{2} \right)   \, - \,  \left( \frac{k-2}{k^2+k+2} \right) \left( r + \ell(k+1)  \right) } \\ \vspace{0.2cm}
& = & \displaystyle{  \frac{1}{2(k^2+k+2)} \left( rk(k+2) - r(k+2) \right)  } \\ \vspace{0.2cm}
&  & \displaystyle{  \hspace*{0.5cm}  \, + \,  \frac{\ell}{2(k^2+k+2)} \left( (k+2)(k^2 + k - 2) - 2(k-2)(k+1)  \right) } \\ \vspace{0.2cm}
& = & \displaystyle{
\frac{2r(k^2+k+2)}{2(k^2+k+2)} \,+ \, \frac{\ell(k^3 + k^2 + 2k)}{2(k^2+k+2)}
 } \\ \vspace{0.2cm}
& = & \displaystyle{ r +  \frac{1}{2} \ell k} \\
& = & \alpha'(G).
\end{array}
\]


Thus, the lower bound in Corollary~\ref{corKeven2} is tight for the class of graphs in the family~$\cG_{k,r}'$.
We show next that Corollary~\ref{corKeven2} is tight for an infinite family of $k$-regular graphs. For $k \ge 4$ an even integer and $r \ge 1$ arbitrary, let $F_{k,r}$ be a graph of order~$n$ and size~$m$ obtained of the disjoint union of $kr/2$ copies of $K_{k+1} - e$ by adding a set $X$ of $r$ new vertices, and adding edges between $X$ and the $kr$ link vertices of degree~$k-1$ in the copies of $K_{k+1} - e$ in such a way that $F_{k,r}$ is a connected, $k$-regular graph. We note that $n = r + \frac{1}{2} k r(k+1) $, $m = rk + \frac{1}{2} k r \times \frac{1}{2}(k^2 + k - 2) = rk + \frac{1}{4} k (k^2 + k - 2) r$ and $\alpha'(G) = r + \frac{1}{2}k r \times \frac{1}{2}k = r + \frac{1}{4}k^2 r$. Therefore the following holds.

\[
\begin{array}{rcl}
\multicolumn{3}{l}{\vspace{0.2cm} \displaystyle{ \left( \frac{k+2}{k^2+k+2} \right) m  \, - \,  \left( \frac{k-2}{k^2+k+2} \right) n  }} \\ \vspace{0.2cm}
& = & \displaystyle{ \left( \frac{r(k+2)}{k^2+k+2} \right) \left( k + \frac{k(k^2 + k - 2)}{4} \right)   \, - \,  \left( \frac{r(k-2)}{k^2+k+2} \right) \left( 1 + \frac{k (k+1)}{2}   \right) } \\ \vspace{0.2cm}
& = & \displaystyle{ \left( \frac{r(k+2)}{k^2+k+2} \right) \left( \frac{k(k^2+k+2)}{4} \right)   \, - \,  \left( \frac{r(k-2)}{k^2+k+2} \right) \left( \frac{k^2 + k + 2}{2} \right) } \\
\vspace{0.2cm}
& = & \displaystyle{ \frac{1}{4}rk(k+2)    \, - \,  \frac{1}{2}r(k-2) } \\
\vspace{0.2cm}
& = & \displaystyle{ r + \frac{1}{4}k^2 r} \\
& = & \alpha'(G).
\end{array}
\]

Thus, Corollary~\ref{corKeven2} is tight for an infinite family of $k$-regular graphs. We state these results formally as follows.

\begin{prop}
The lower bound in Corollary~\ref{corKeven2} is achieved for both $k$-regular graphs and for the class of graphs in the family~$\cG_{k,r}'$. \label{p:Gkrregular}
\end{prop}

\section{The convex set $L_k$}

Let $k \ge 3$ be an integer. Let $\cG_k$ denote the class of connected graphs with maximum degree at most~$k$. For every pair $(a,b)$ of real numbers $a$ and
$b$ we define the concept of $k$-good, $k$-bad and $k$-tight as follows.

\begin{itemize}
 \item $(a,b)$ is called $k$-\emph{good} if there exists a constant $T_{a,b}$ such that
\[
\alpha'(G) \ge a |V(G)| + b |E(G)| - T_{a,b}
\]
holds for all $G \in \cG_k$.

 \item $(a,b)$ is called $k$-\emph{bad} if it is not $k$-good.

 \item $(a,b)$ is called $k$-\emph{tight} if it is $k$-good and there exists a constant $S_{a,b}$ such that
\[
\alpha'(G) \le a |V(G)| + b |E(G)| - S_{a,b}
\]
holds for  infinitely many graphs $G \in \cG_k$.
\end{itemize}

If we say that $(a,b)$ is $k$-tight for a certain subset of $\cG_k$ (for example, the class of trees or $k$-regular graphs), then we mean that there are infinitely many graphs
from this class that satisfy $\alpha'(G) \le a |V(G)| + b |E(G)| - S_{a,b}$ for some constant $S_{a,b}$.

Suppose that $(a,b)$ is $k$-good and $\varepsilon \ge 0$. Then, there exists a constant $T_{a,b}$ such that $\alpha'(G) \ge a |V(G)| + b |E(G)| - T_{a,b} \ge a |V(G)| + (b - \varepsilon) |E(G)| - T_{a,b}$, implying that $(a,b - \varepsilon)$ is $k$-good where $T_{a,b- \varepsilon} = T_{a,b}$. We state this formally as follows.

\begin{ob}
If $(a,b)$ is $k$-good and $\varepsilon \ge 0$, then $(a,b-\varepsilon)$ is $k$-good.
\label{ob1}
\end{ob}

\begin{lem}\label{lem_pic1}
If $(a,b)$ is $k$-good and $\varepsilon \ge 0$, then both $(a+ \varepsilon,b-\varepsilon)$ and $(a- \varepsilon \cdot k, b +  2 \varepsilon)$ are $k$-good.
Furthermore the following holds.
\\[-28pt]
\begin{enumerate}
\item If $(a,b)$ is $k$-tight for trees, then $(a+ \varepsilon,b-\varepsilon)$ is $k$-tight.
\item If $(a,b)$ is $k$-tight for $k$-regular graphs, then $(a - \varepsilon \cdot k, b + 2 \varepsilon)$ is $k$-tight.
\end{enumerate}
\end{lem}
\proof
Let $G \in \cG_k$ have order $n$ and size $m$.  Since $G$ is connected, we note that $m \ge n-1$. Since $(a,b)$ is $k$-good, this implies that there exists a constant $T_{a,b}$ such that the following also holds for $\varepsilon \ge 0$.
\[
\begin{array}{rcl} \vspace{0.2cm}
\alpha'(G) & \ge & a \cdot n + b \cdot m - T_{a,b} \\ \vspace{0.2cm}
  & \ge & a \cdot n + b \cdot m - T_{a,b} + \varepsilon (n-1 - m) \\ \vspace{0.2cm}
  & = & (a + \varepsilon) n  + (b-\varepsilon) m - (T_{a,b} + \varepsilon).
\end{array}
\]

So letting $T_{a+ \varepsilon,b-\varepsilon} = T_{a,b} + \varepsilon$, the pair $(a+ \varepsilon,b-\varepsilon)$ is $k$-good.
If $(a,b)$ is $k$-tight for trees, then there exists a constant $S_{a,b}$ such that for infinitely many trees, $G'$, in $\cG_k$ we have $\alpha'(G') \le a |V(G')| + b |E(G')| - S_{a,b}$. Let $G'$ has order~$n'$ and size~$m'$. Then, $m' = n' - 1$ and, analogously as before, the following holds.
\[
\alpha'(G') \le a \cdot n' + b \cdot m' - S_{a,b} = (a + \varepsilon) n'  + (b-\varepsilon) m' - (S_{a,b} + \varepsilon).
\]

So letting $S_{a+ \varepsilon,b-\varepsilon} = S_{a,b} + \varepsilon$,  the pair $(a+ \varepsilon,b-\varepsilon)$ is $k$-tight in this case.
Recall that $G\in \cG_k$ has order $n$ and size $m$. As $G$ has maximum degree at most~$k$, we have $nk \ge \sum_{v \in V(G)} d_G(v) = 2m$, which implies
that the following also holds  for all $\varepsilon \ge 0$.
\[
\begin{array}{rcl} \vspace{0.2cm}
\alpha'(G) & \ge & a \cdot n + b \cdot m - T_{a,b} \\ \vspace{0.2cm}
  & \ge & a \cdot n + b \cdot m - T_{a,b} + \varepsilon (2m-nk) \\ \vspace{0.2cm}
  & = & (a - \varepsilon \cdot k) n  + (b + 2 \varepsilon) m - T_{a,b}.
\end{array}
\]

So letting $T_{a - \varepsilon \cdot k, b + 2 \varepsilon} = T_{a,b}$, we note that $(a - \varepsilon \cdot k, b + 2 \varepsilon)$ is $k$-good.
If $(a,b)$ is $k$-tight for $k$-regular graphs, then for infinitely many $k$-regular graphs, $G'$, in $\cG_k$ we have
$\alpha'(G') \le a |V(G')| + b |E(G')| - S_{a,b}$. Let $G'$ has order~$n'$ and size~$m'$. Then, $n'k = 2m'$ and, analogously as before, the following holds.
\[
\begin{array}{rcl} \vspace{0.2cm}
\alpha'(G) & \le & a \cdot n' + b \cdot m' - S_{a,b} \\ \vspace{0.2cm}
  & = & a \cdot n' + b \cdot m' - S_{a,b} + \varepsilon (2m'-n'k) \\ \vspace{0.2cm}
  & = & (a - \varepsilon \cdot k) n'  + (b + 2 \varepsilon) m' - S_{a,b}.
\end{array}
\]
So letting $S_{a- \varepsilon \cdot k, b + 2 \varepsilon} = S_{a,b}$,  the pair $(a - \varepsilon \cdot k,b + 2 \varepsilon)$ is $k$-tight in this case.~\qed

\begin{lem}\label{lem_pic2}
If $(a,b)$ is $k$-tight, then $(a+\varepsilon,b)$ and $(a,b+\varepsilon)$ are both $k$-bad for all $\varepsilon > 0$.
\end{lem}
\proof
Assume that $(a,b)$ is $k$-tight and, for the sake of contradiction, suppose that $(a+\varepsilon,b)$ is $k$-good for some $\varepsilon > 0$.
That is, there exists a constant $T_{a+\varepsilon,b}$ such that
\[
\alpha'(G)  \ge  (a+\varepsilon) \cdot |V(G)| + b \cdot |E(G)| - T_{a+\varepsilon,b}
\]

\noindent holds for all $G \in \cG_k$.
Since $(a,b)$ is $k$-tight, there exists a constant $S_{a,b}$ such that
\[
\alpha'(G)   \le  a \cdot |V(G)| + b \cdot |E(G)| - S_{a,b}
   =  (a + \epsilon) \cdot |V(G)| + b \cdot |E(G)| - S_{a,b} - \varepsilon |V(G)|
\]

\noindent
holds for infinitely many $G \in \cG_k$. However as there are infinitely many such graphs, we can choose such a graph $G$ (of sufficiently large order) such that $S_{a,b} + \varepsilon |V(G)| > T_{a+\varepsilon,b}$.
For this graph $G$, we have $\alpha'(G) < (a + \varepsilon) |V(G)| + b |E(G)| - T_{a+\varepsilon,b}$, a contradiction. Therefore, $(a+\varepsilon,b)$ is $k$-bad.

The fact that $(a,b+\varepsilon)$ is $k$-bad can be proved analogously.~\qed

\begin{lem}\label{lem_pic3}
If $(a_1,b_1)$ and $(a_2,b_2)$ are both $k$-good,  then
$(\varepsilon a_1 + (1-\varepsilon) a_2, \varepsilon b_1 + (1-\varepsilon) b_2)$ is also $k$-good for all $0 \le \varepsilon \le 1$.

Furthermore if $(a_1,b_1)$ and $(a_2,b_2)$ are both $k$-tight for the same infinite class ${\cal G'} \subseteq \cG_k$, then
$(\varepsilon a_1 + (1-\varepsilon) a_2, \varepsilon b_1 + (1-\varepsilon) b_2)$ is also $k$-tight.
\end{lem}
\proof
Since $(a_1,b_1)$ and $(a_2,b_2)$ are both $k$-good, there exists constants $T_{a_1,b_1}$ and $T_{a_2,b_2}$ such that
\[
\alpha'(G) \ge a_1 |V(G)| + b_1 |E(G)| - T_{a_1,b_1} \hspace*{0.5cm}  \mbox{and} \hspace*{0.5cm} \alpha'(G) \ge a_2 |V(G)| + b_2 |E(G)| - T_{a_2,b_2}
\]

\noindent hold for all $G \in \cG_k$. Let $0 \le \varepsilon \le 1$. Multiplying the first equation by $\varepsilon$ and the second by $(1-\varepsilon)$, and then adding the equations together shows that
\[
\alpha'(G) \ge (\varepsilon a_1 + (1-\varepsilon) a_2) |V(G)| + (\varepsilon b_1 + (1-\varepsilon) b_2) |E(G)|
- \varepsilon T_{a_1,b_1} - (1-\varepsilon) T_{a_2,b_2}
\]

\noindent holds for all $G \in \cG_k$. This proves that $(\varepsilon a_1 + (1-\varepsilon) a_2, \varepsilon b_1 + (1-\varepsilon) b_2)$ is $k$-good, with
$T_{\varepsilon a_1 + (1-\varepsilon) a_2, \varepsilon b_1 + (1-\varepsilon) b_2} = \varepsilon T_{a_1,b_1} + (1-\varepsilon) T_{a_2,b_2}$.

Assume that $(a_1,b_1)$ and $(a_2,b_2)$ are both $k$-tight for the same infinite class ${\cal G'} \subseteq \cG_k$. Thus there exists constants $S_{a_1,b_1}$ and $S_{a_2,b_2}$ such that
\[
\alpha'(G) \le a_1 |V(G)| + b_1 |E(G)| - S_{a_1,b_1} \hspace*{0.5cm}  \mbox{and} \hspace*{0.5cm} \alpha'(G) \le a_2 |V(G)| + b_2 |E(G)| - S_{a_2,b_2}
\]

\noindent hold for all $G \in \cG'$. Again, multiplying the first equation by $\varepsilon$ and the second by $(1-\varepsilon)$ and adding the equations together shows that
\[
\alpha'(G) \le (\varepsilon a_1 + (1-\varepsilon) a_2) |V(G)| + (\varepsilon b_1 + (1-\varepsilon) b_2) |E(G)|
- \varepsilon S_{a_1,b_1} - (1-\varepsilon) S_{a_2,b_2}.
\]

\noindent holds for all $G \in \cG'$. Therefore, $(\varepsilon a_1 + (1-\varepsilon) a_2, \varepsilon b_1 + (1-\varepsilon) b_2)$ is $k$-tight letting
\[
S_{\varepsilon a_1 + (1-\varepsilon) a_2, \varepsilon b_1 + (1-\varepsilon) b_2} = \varepsilon S_{a_1,b_1} + (1-\varepsilon) S_{a_2,b_2}. \hspace*{0.5cm} \Box
\]

\subsection{$k$ odd}

\begin{thm} \label{thm_pic_odd}
Let $k \ge 3$ be odd and let $a^* = \frac{k-1}{k(k^2 - 3)}$ and $b^* = \frac{k^2 - k - 2}{k(k^2 - 3)}$. For any pair $(a,b)$, the following holds.
\\[-28pt]
\begin{enumerate}
\item If $a \le a^*$, then $(a,b)$ is $k$-good if and only if $b \le b^*+ 2 \left( \frac{a^*-a}{k} \right)$.
 \item If $a > a^*$, then $(a,b)$ is $k$-good if and only if $b \le b^*+ a^*-a$.
\end{enumerate}
\end{thm}
\proof
By Corollary~\ref{corKodd}, the pair $(a^*,b^*)$ is $k$-good with
\[
T_{a^*,b^*} = \frac{k-1}{k(k^2 - 3)}.
\]

By Proposition~\ref{p:Hkrtreeregular}, the lower bound in Corollary~\ref{corKodd} is achieved for an infinite class of both trees and $k$-regular graphs, implying that $(a^*,b^*)$ is $k$-tight for both trees and $k$-regular graphs.

Suppose $a \le a^*$, and let $\varepsilon = (a^*-a)/k$. By Lemma~\ref{lem_pic1}, we note that $(a^*-\varepsilon \cdot k, b^* + 2 \varepsilon)$ is $k$-good. Further, since $(a^*,b^*)$ is $k$-tight for $k$-regular graphs, by Lemma~\ref{lem_pic1}(b), we note that $(a^*-\varepsilon \cdot k, b^* + 2 \varepsilon)$ is $k$-tight. Since $\varepsilon = (a^*-a)/k$, this is equivalent to $(a,b^*+ 2 \left( \frac{a^*-a}{k} \right))$ being $k$-tight. If $b \le  b^* + 2 \left( \frac{a^*-a}{k} \right)$, then by Observation~\ref{ob1}, the pair $(a,b)$ is $k$-good. If $b > b^* + 2 \left( \frac{a^*-a}{k} \right)$, then by Lemma~\ref{lem_pic2} and our earlier observation that $(a,b^*+ 2 \left( \frac{a^*-a}{k} \right))$ is $k$-tight, the pair $(a,b)$ is $k$-bad. This completes the case when $a \le a^*$.

Suppose next that $a > a^*$, and let $\varepsilon = a-a^*$. By Lemma~\ref{lem_pic1}, we note that $(a^* + \varepsilon, b^* - \varepsilon)$ is $k$-good. Further since $(a^*,b^*)$ is $k$-tight for trees, we note by Lemma~\ref{lem_pic1}(a) that $(a^* + \varepsilon, b^* - \varepsilon)$ is $k$-tight. Since $\varepsilon = (a^*-a)/k$, this is equivalent to $(a,b^* + a^*-a)$ being $k$-tight. If $b \le  b^* + a^*-a$, then by Observation~\ref{ob1}, the pair $(a,b)$ is $k$-good. If $b > b^* + a^*-a$, then by Lemma~\ref{lem_pic2} and our earlier observation that $(a,b^* + a^*-a)$ is $k$-tight, the pair $(a,b)$ is $k$-bad.~\qed

\medskip
We remark that the equation in Theorem~\ref{thm_pic_odd}(a) corresponds to the half-plane $\ell_2$ described in the introductory section, noting that
\[
b \le b^*+ 2 \left( \frac{a^*-a}{k} \right) = - \left(\frac{2}{k} \right) a \, +  \, \frac{k^3 - k^2 - 2}{k^2(k^2 - 3)}.
\]

\noindent
The equation in Theorem~\ref{thm_pic_odd}(b) corresponds to the half-plane $\ell_1$ described in the introductory section, noting that
\[
b \le b^*+ a^*-a  = - a \, +  \, \frac{1}{k}. 
\]

Theorem~\ref{thm_pic_odd} is illustrated in Figure~\ref{fig:Lk} when $k=3$ and $k=5$. The grey area corresponds to all $k$-good pairs $(a,b)$ while the non-grey area
corresponds to the $k$-bad pairs.

\newpage
\subsection{$k$ even}

\begin{thm} \label{thm_pic_even}
Let $k \ge 4$ be even and let
$a_1^* = \frac{1}{k(k+1)}$  and $b_1^* = \frac{1}{k+1}$ and $a_2^* = - \frac{k-2}{k^2+k+2}$  and $b_2^* = \frac{k+2}{k^2+k+2}$. For any pair $(a,b)$, the following holds.
\\[-22pt]
\begin{enumerate}
\item If $a \le a_2^*$, then $(a,b)$ is $k$-good if and only if $b \le  b_2^* + 2\left(\frac{a^*-a}{k} \right)$.
\item If $a > a_1^*$, then $(a,b)$ is $k$-good if and only if $b \le b_1^* + a_1^*-a$.
\item If $a_2^* < a \le a_1^*$, then $(a,b)$ is $k$-good if and only if $b \le b_2^* + \frac{(b_1^*  - b_2^*)(a-a_2^*)}{a_1^* - a_2^*} $.
\end{enumerate}
\end{thm}
\proof
By Corollary~\ref{corKeven}, the pair $(a_1^*,b_1^*)$ is $k$-good with
$T_{a_1^*,b_1^*} = - 1/k$, while by Corollary~\ref{corKeven2}, the pair $(a_2^*,b_2^*)$ is $k$-good with
\[
T_{a_2^*,b_2^*} = - \frac{k+2}{k^2 + k + 2}.
\]

By Proposition~\ref{p:Gkrtree}, the lower bound in Corollary~\ref{corKeven} is achieved for both trees and for the class of graphs in the family~$\cG_{k,r}'$, implying that $(a_1^*,b_1^*)$ is $k$-tight for both trees and graphs in the family~$\cG_{k,r}'$.
By Proposition~\ref{p:Gkrregular}, the lower bound in Corollary~\ref{corKeven2} is achieved for the class of graphs in the family~$\cG_{k,r}'$ and for the class of $k$-regular graphs, implying that $(a_2^*,b_2^*)$ is $k$-tight for these classes of graphs.

Suppose that $a \le a_2^*$, and let $\varepsilon = (a_2^*-a)/k$. By Lemma~\ref{lem_pic1}, we note that $(a_2^*-\varepsilon \cdot k, b_2^* + 2 \varepsilon)$ is $k$-good. Further, since $(a_2^*,b_2^*)$ is $k$-tight for $k$-regular graphs, by Lemma~\ref{lem_pic1}(b), we note that $(a_2^*-\varepsilon \cdot k, b_2^* + 2 \varepsilon)$ is $k$-tight. Since $\varepsilon = (a_2^*-a)/k$, this is equivalent to
$(a,b_2^*+ 2 \left( \frac{a_2^*-a}{k} \right))$ being $k$-tight.
If $b \le  b_2^* + 2 \left( \frac{a^*-a}{k} \right)$, then by Observation~\ref{ob1}, the pair $(a,b)$ is $k$-good.
If $b > b_2^* + 2 \left( \frac{a^*-a}{k} \right)$, then by Lemma~\ref{lem_pic2} and our earlier observation that $(a,b_2^*+ 2 \left( \frac{a_2^*-a}{k} \right))$ is $k$-tight, the pair $(a,b)$ is $k$-bad.  This completes the case when $a \le a_2^*$.

Suppose that $a > a_1^*$, and let $\varepsilon = a_1-a^*$. By Lemma~\ref{lem_pic1}, we note that $(a_1^* + \varepsilon, b_1^* - \varepsilon)$ is $k$-good. Further since $(a_1^*,b_1^*)$ is $k$-tight for trees, we note by Lemma~\ref{lem_pic1}(a) that $(a_1^* + \varepsilon, b_1^* - \varepsilon)$ is $k$-tight. Since $\varepsilon = (a^*-a)/k$, this is equivalent to $(a,b_1^* + a_1^*-a)$ being $k$-tight.
If $b \le  b_1^* + a_1^*-a$, then by Observation~\ref{ob1}, the pair $(a,b)$ is $k$-good.
If $b > b_1^* + a_1^*-a$, then by Lemma~\ref{lem_pic2} and our earlier observation that $(a,b_1^* + a_1^*-a)$ is $k$-tight, the pair $(a,b)$ is $k$-bad.

Finally, suppose that $a_2^* < a \le a_1^*$ and let $\varepsilon = (a-a_2^*)/(a_1^* - a_2^*)$. By Lemma~\ref{lem_pic3}, we note that
$(\varepsilon a_1^* + (1- \varepsilon) a_2^*, \varepsilon b_1^* + (1-\varepsilon) b_2^*)$ is $k$-good. Furthermore, since $(a_1^*,b_1^*)$ and $(a_2^*,b_2^*)$ are both $k$-tight for graphs in the family~$\cG_{k,r}'$, we note that $(\varepsilon a_1^* + (1- \varepsilon) a_2^*, \varepsilon b_1^* + (1-\varepsilon) b_2^*)$ is $k$-tight. Since $\varepsilon = (a-a_2^*)/(a_1^* - a_2^*)$, this is equivalent to $\left( a, b_2^* + \frac{(b_1^*  - b_2^*)(a-a_2^*)}{a_1^* - a_2^*} \right)$ being $k$-tight.
If $b \le  b_2^* + \frac{(b_1^*  - b_2^*)(a-a_2^*)}{a_1^* - a_2^*}$, then by Observation~\ref{ob1}, the pair $(a,b)$ is $k$-good.
If $b > b_2^* + \frac{(b_1^*  - b_2^*)(a-a_2^*)}{a_1^* - a_2^*}$, then by Lemma~\ref{lem_pic2} and our earlier observation that $\left( a, b_2^* + \frac{(b_1^*  - b_2^*)(a-a_2^*)}{a_1^* - a_2^*} \right)$ is $k$-tight, the pair $(a,b)$ is $k$-bad.~\qed

\medskip
We remark that the equation in Theorem~\ref{thm_pic_even}(a) corresponds to the half-plane $\ell_3$ described in the introductory section, noting that
\[
b \le b_2^* + 2\left(\frac{a_2^*-a}{k} \right) = - \left( \frac{2}{k} \right) a \, +  \, \frac{k^2 + 4}{k(k^2 + k + 2)}.
\]

\noindent 
The equation in Theorem~\ref{thm_pic_even}(b) corresponds to the half-plane $\ell_1$ described in the introductory section, noting that
\[
b \le b_1^* + a_1^*-a  = - a \, +  \, \frac{1}{k}.
\]

\noindent 
The equation in Theorem~\ref{thm_pic_even}(c) corresponds to the half-plane $\ell_4$ described in the introductory section, noting that
\[
b \le b_2^* + \frac{(b_1^*  - b_2^*)(a-a_2^*)}{a_1^* - a_2^*} = -\left(\frac{2k^2 }{k^3-k+2} \right) a \, +  \, \frac{k^2 - k + 2}{k^3-k+2}.
\]

Theorem~\ref{thm_pic_even} is illustrated in Figure~\ref{fig:Lk} when $k=4$ and $k=6$. The grey area corresponds to all $k$-good pairs $(a,b)$ while the non-grey area corresponds to the $k$-bad pairs.


\medskip
\begin{thebibliography}{99}

\bibitem{Berge} C. Berge, \textit{C. R. Acad. Sci. Paris Ser. I Math.} \textbf{247,} (1958) 258--259 and \textit{Graphs and Hypergraphs} (Chap. 8, Theorem 12), North-Holland, Amsterdam, 1973.

\bibitem{Biedl} T. Biedl, E. D. Demaine, C. A. Duncan, R. Fleischer and S. G. Kobourov, Tight bounds on maximal and maximum matchings. \textit{Discrete Math.} \textbf{285} (2004), 7--15.

\bibitem{ChMc} V. Chv\'{a}tal and C. McDiarmid, Small transversals in hypergraphs. \textit{Combinatorica} \textbf{12} (1992), 19--26.

\bibitem{CiGrHa09} S. M. Cioab\u{a}, D. A. Gregory, and W. H. Haemers, Matchings in regular graphs from eigenvalues. \textit{J. Combinatorial Theory Ser. B} \textbf{99} (2009), 287--297.

\bibitem{HaSc14} P. E. Haxell and A. D. Scott, On lower bounds for the matching number of subcubic graphs. Manuscript, June 2014. http://arxiv.org/abs/1406.7227


\bibitem{HeLoRa12} M. A. Henning, C. L\"{o}wenstein, and D. Rautenbach, Independent sets and matchings in subcubic graphs. \textit{Discrete Math.} \textbf{312} (2012), 1900--1910.

\bibitem{HeYe07} M. A. Henning and A. Yeo, Tight lower bounds on the size of a matching in a regular graph. \textit{Graphs Combin.} \textbf{23} (2007), 647--657.


\bibitem{MHAYbookTD} M. A. Henning and A. Yeo, \emph{Total domination in graphs (Springer Monographs in Mathematics)}.  ISBN-13: 978-1461465249  (2013).

\bibitem{JaWe13} S. Jahanbekam and D. B. West, New lower bounds on matching numbers of general and bipartite graphs. \textit{Congr. Numer.} \textbf{218} (2013), 57--–59.



\bibitem{Ka87} P. Katerinis, Maximum matching in a regular graph of specified connectivity and bounded order. \textit{J. Graph Theory} \textbf{11} (1987), 53--58.

\bibitem{LoPl86} L. Lov\'{a}sz and M. D. Plummer, Matching Theory, North-Holland Mathematics Studies, vol. 121, Ann. Discrete Math., vol. 29, North-Holland, 1986.


\bibitem{OWest10} Suil O and D. B. West, Balloons, cut-edges, matchings, and total domination in regular graphs of odd degree. \textit{J. Graph Theory} \textbf{64} (2010), 116--131.


\bibitem{OWest11} Suil O and D. B. West, Matching and edge-connectivity in regular graphs. \textit{European J. Comb.} \textbf{32} (2011), 324--329.

\bibitem{Pl03} M. Plummer, Factors and Factorization. 403--430. {\em Handbook of Graph Theory} {\em ed.} J. L. Gross and J. Yellen. CRC Press, 2003, ISBN: 1-58488-092-2.

\bibitem{Pu95} W. R. Pulleyblank, Matchings and Extension. 179--232. {\em Handbook of Combinatorics} {\em ed.} R. L. Graham, M. Gr\"{o}tschel, L. Lov\'{a}sz. Elsevier Science B.V. 1995, ISBN 0-444-82346-8.


\bibitem{West11} D. B. West, A short proof of the Berge-Tutte Formula and the Gallai-Edmonds Structure Theorem. \textit{European J. Comb.} \textbf{32} (2011), 674--676.

\end{thebibliography}
\end{document}